\documentclass[12pt]{article}

\usepackage{amssymb}

\usepackage{latexsym,amsmath,amsfonts,amscd}
\usepackage{xfrac}
\usepackage{graphics}
\usepackage{caption}
\usepackage{subcaption}
\usepackage{epstopdf}
\usepackage{color}
\usepackage{changebar}
\usepackage{indentfirst}
\usepackage{verbatim}
\usepackage{xcolor}
\usepackage{hyperref}
\usepackage{amssymb}
\usepackage{lscape}
\usepackage{array,multirow}
\usepackage{booktabs}

\newcommand{\first}[1]{\textcolor{black}{#1}}
\newcommand{\second}[1]{\textcolor{black}{#1}}

\newtheorem{remark}{Remark}

\usepackage{authblk}
\usepackage[margin=3.0cm]{geometry}

\title{A PINN approach for the online identification and control of unknown PDEs }

\author[$\dagger$]{Alessandro Alla \footnote{Corresponding author. Email address: \textit{alessandro.alla@uniroma1.it}}}
\author[$\ddagger$]{Giulia Bertaglia}
\author[$\mathsection$]{Elisa Calzola}
\affil[$\dagger$]{\small Department of Mathematics,
Sapienza University of Rome, Italy}
\affil[$\ddagger$]{\small Department of Environmental and Prevention Sciences,
University of Ferrara, Italy}
\affil[$\mathsection$]{\small Department of Computer Science,
University of Verona, Italy}

\begin{document}

\maketitle

\begin{abstract}
Physics-Informed Neural Networks (PINNs) have revolutionized solving differential equations by integrating physical laws into neural networks training. This paper explores PINNs for open-loop optimal control problems (OCPs) with incomplete information, such as \first{sparse} initial and boundary data \first{and partially unknown system parameters}. We derive optimality conditions from the Lagrangian multipliers and use PINNs to predict the state, adjoint, and control variables. In contrast with previous methods, our approach integrates these elements into a single neural network and addresses scenarios with consistently limited data. In addition, we address the study of partially unknown equations identifying underlying parameters online by searching for the optimal solution \first{recurring to a 2-in-series architecture of PINNs, in which scattered data of the uncontrolled solution is used}. Numerical examples show the effectiveness of the proposed method even in scenarios characterized by a considerable lack of information.
\end{abstract}


\tableofcontents

\section{Introduction}
\label{sec:intro}

Physics-Informed Neural Networks (PINNs) have recently revolutionized the use of neural networks for solving differential equations \cite{E2021,E2020}. Unlike traditional Deep Neural Networks (DNNs), which rely solely on user-provided data for the training, PINNs incorporate additional loss terms that encode the residuals of the differential equations underlying the physical phenomenon under investigation. This extension allows neural networks to reconstruct complex dynamical systems more consistently \cite{Baker2019}. Moreover, when set up knowledgeably, the PINN approach has shown robust enough to also allow reliable predictions of evolutionary dynamics over time, even when dealing with complex multiscale systems \cite{Bertaglia2022b,Bertaglia2022a,Jin2021,Jin2022}. On the other hand, PINNs have also been used to solve inverse problems by identifying the differential equations governing the dynamics of interest from the available data. Indeed, the contribution of PINNs to the scientific community is twofold: (i) approximating forward problems and (ii) solving inverse problems.
The effectiveness of this approach has been demonstrated in numerous studies. For foundational works on this topic, we refer to \cite{Karniadakis2021,RPK19}.

Among others, optimal control problems (OCPs) can be considered a special case of inverse problems, where for a given differential equation, a space-time function $u(x,t)$ is unknown and needs to be computed or discovered to satisfy a desired configuration. In the mathematical community, OCPs for partial differential equations (PDEs) have a solid mathematical foundation with strong theoretical results. In fact, many efficient numerical methods have been developed to solve this kind of problem (we refer e.g. to \cite{HPUU09} and references therein). However, traditional numerical approaches are based on the assumption that initial conditions, boundary conditions, and modeling parameters are all precisely known, which is frequently not the case, especially for real-world applications.

Inspired by the particular structure of neural networks, which by their nature are trained following a dual flow of forward (feed-forward) and backward (back-propagation) types, reminiscent precisely of the combined flow typical of OCPs, in this paper we explore the usage of PINNs for solving optimal control problems.
To the best of our knowledge, \cite{MN23} is the first work to propose using PINNs for such  problems. In that work, the authors make use of different neural networks and the loss function also includes the cost function of the problem, which makes their approach different from what is proposed in the present manuscript. Moreover, in their study, the control is computed as an inverse problem, thus it is treated as a neural network parameter and not as an output variable.

This field remains largely open and dynamic, with most contributions emerging only in the last couple of years. This highlights the novelty and active research interest in this area. Among the most recent works, we mention that in \cite{GKP23} the authors formulate the minimization of a loss function involving the residual of the PDE together with the initial, boundary, and final conditions. They also provide error estimates for the training error, demonstrating that it diminishes as the size of the neural network and the number of training points increase to infinity. 
Parametric OCPs have also been studied in \cite{DSR23,YXTY24}. Specifically, in \cite{YXTY24} an adjoint-oriented neural network method for all-at-once solutions of parametric optimal control is discussed. The authors introduce three different neural networks, each corresponding to one of the equations in the optimality system. Their method employs an iterative approach within a supervised learning framework, using a Monte Carlo strategy for the sampling. Nevertheless, all simulations in this study are conducted on \textit{stationary} PDEs. 

In the present work, we propose to extend the applicability of PINNs to solve time-dependent optimal control problems using an open-loop approach. Specifically, we focus on scenarios where complete information on the problem is not available. For example, the initial and boundary conditions for the state equation can be prescribed only at a sparse level. These circumstances make the application of traditional numerical methods essentially impossible. 
In addition, we consider scenarios in which the state and adjoint equations are partially unknown and the model parameters need to be discovered, thus testing the methodology in solving inverse problems for finding the unknown parameters along with solving the OCP. This provides an \textit{online} algorithm in which we identify and control, at the same time, the partially unknown problem. 
Thus, the proposed method leverages a small amount of data to exploit the features of PINNs in solving both the inverse and forward problems investigated. 
Despite the challenging context, we demonstrate how the proposed OCP-PINNs\footnote{OCP-PINNs stands for Optimal Control Problem using Physics Informed Neural Networks.} can effectively generate informative solutions even in such scenarios, and we also show the capability of the method to extrapolate information beyond the time frame provided by the available data. \first{Although adjoint information is used in the method presented here, as in \cite{DSR23}, the approach proposed here differs from that in \cite{DSR23} because here we are dealing with time-dependent problems and also with the identification of the control variable.} \second{Furthermore, we would like to underline that the methodology here introduced also differs from standard Reinforcement Learning strategies \cite{BS18}, because here we are proposing an online approach in which the observations are given and will not change during the tuning of the neural network.}

Mathematically, we derive the optimality conditions from the method of the Lagrangian multipliers, which leads us to a system of three equations: the state equation for the evolution of the state variable \(y\), the adjoint equation for the evolution of the adjoint variable \(p\), and the optimality condition that relates the control $u$ to the adjoint.

For the sake of completeness, we also mention other methods aimed at identifying differential equations. Sparse Identification of Nonlinear Dynamics (SINDy, \cite{RABK19,RBPK16}) has been applied to both ODEs \cite{BPK16} and PDEs. This method uses sparse optimization techniques to identify the model based on a chosen library of functions, which typically includes derivative terms and nonlinear functions that might appear in the equation. SINDy has also been utilized for identifying controlled systems, as shown in \cite{KKB19}, where an external source is used as input to identify the system, followed by the application of Nonlinear Model Predictive Control (NMPC) for system regulation. Some recent works, \cite{AP24,APPP23,zukb24}, have extended the SINDy strategy for online identification and control of PDEs. Additionally, other strategies focus on simultaneous control and identification of systems, such as the approach proposed in \cite{KS08} for PDEs and in \cite{MLG20} for ODEs. A recent study in \cite{CFDKMRV21} explores controlling unknown systems using Model Predictive Control (MPC), where system identification is performed using Extended Dynamic Mode Decomposition (EDMD), representing a surrogate linear model.

The rest of the paper is organized as follows. In Section \ref{sec:NN-PINN}, we provide a summary of the basics for constructing a PINN. Section \ref{sec:PINN-OCP} is devoted to presenting the here proposed OCP-PINNs approach, hence the extension of the PINN framework for solving OCPs using Lagrange multipliers. Several numerical tests are presented in Section \ref{sec:tests} to demonstrate the validity of the methodology. Final discussion and conclusions are given in Section \ref{sec:final}.

\section{The PINN framework}
\label{sec:NN-PINN}
Let us assume we want to evaluate the space-time dynamics  $y(x,t)\in \Omega =  Q\times [0,T] \subset \mathbb{R}^{d_x}\times\mathbb{R}$, through a data-driven approach, recurring to deep neural networks. For example, we consider an $(L+1)$ layered feed-forward neural network which consists of an input layer, an output layer, and $(L-1)$ hidden layers, and we focus on the 1-D case in space, being $d_x=1$. Given an available dataset of $N_d$ observations located in the spatio-temporal domain respectively in $(x_d^i,t_d^i) \subset \Omega$, $i=1,\ldots N_d$ as input, the DNN can be used to output a prediction of the value $  y_{NN}(x,t; \theta)\approx y(x,t)$ parameterized by the set of neural network parameters $\theta$. To this aim, if we indicate with $z=(x,t)$ the input vector, we can define the DNN as follows:
\begin{eqnarray*}
z^1 &=& W^1 z + b^1\,,\\
z^l &=& \sigma \circ (W^l z^{l-1} + b^l),\quad l=2,\ldots,L-1\,,\\
y_{NN}(z;\theta) = z^L&=& W^L z^{L-1} + b^L\,,
\end{eqnarray*}
where $W^{l}\in\mathbb{R}^{n_{l}\times n_{l-1}}$ is the matrix of the weights of the neural network, while $b^l\in\mathbb{R}^{n_{l}}$ is the vector of the biases, which form together the set of parameters that the DNN needs to ``learn": $\theta=(W^1,b^1,\ldots,W^L,b^L).$ We refer to e.g. \cite{goodfellow2016} for a more detailed description. Here, \second{$n_l$} is the width of the $l$-th hidden layer, with $n_1$ the input dimension and $n_L$ the output dimension,  while $\sigma$ is the activation function (such as the ReLU function, the \second{sigmoid} function, the \texttt{tanh}, ...).

To learn the model $y_{NN}$, the network’s free parameters $\theta$ are tuned to find the optimal values at each iteration of the training phase through a supervised learning process so that the DNN’s predictions closely match the available experimental data. This is usually done by minimizing a loss function $\mathcal{L}_d(\theta)$ (also called cost or risk function) that, in the simplest case of standard DNNs, consists in the sole mean-squared-error (MSE) between the DNN's model evaluated in the training points $(x_d,t_d)$ and the real dataset. For our setting, the loss term evaluating the mismatch to the dataset reads:
\begin{equation}
\label{eq:general-nn-loss}
\mathcal{L}_d(\theta) = \frac{1}{N_d} \sum_{i=1}^{N_d} \left| y_{NN}(x_d^i, t_d^i;\theta) - y(x_d^i, t_d^i) \right|^2\,,
\end{equation}
where $\{y(x_d^i, t_d^i)\}_{i=1}^{N_d}$ is the available measurement dataset (and the subscript ``$d$" indicates the dependence on the data).
Notice that, in practice, the performance of the neural network is estimated on a finite dataset (which is unrelated to any data used to train the model) and called \textit{test error} (or validation error), whereas the error in the loss function (which is used for training purposes) is called the \textit{training error}.

As already discussed in Section \ref{sec:intro}, using a purely data-driven approach in which the DNN is trained only through $\mathcal{L}_d$ can have significant weaknesses. Indeed, whilst a standard DNN might accurately model the physical process within the vicinity of the experimental data, it will normally fail to generalize away from the training points $(x_d^i, t_d^i)$. This happens because the neural network has not \textit{truly} learned the physical dynamics of the phenomenon \cite{Karniadakis2021}.

To improve the effectiveness of the learning process, PINNs take advantage of the prior physical/scientific knowledge we possess including it in the neural network workflow. 
More specifically, let us assume that, in addition to having access to experimental data, we know that the dynamic investigated is governed by a differential operator $\mathcal{F}(y, x, t) = 0$ in the spatio-temporal domain $\Omega \subset \mathbb{R}\times\mathbb{R}$. Furthermore, we may know that the dynamic is subject to a general operator $\mathcal{B}(y, x, t) = 0$ that prescribes arbitrary initial and boundary conditions of the system in $\partial \Omega$, at specific (sparse) locations. 

To transform a standard DNN into a PINN, it is sufficient to require the neural network to satisfy not only the agreement with experimental data (as for \eqref{eq:general-nn-loss}) but also the physical laws we know of the system, acting directly in the loss function $\mathcal{L}(\theta)$ by including additional terms: 
\begin{equation}
\label{eq:general-pinn-loss}
\mathcal{L}(\theta) = w_d \mathcal{L}_{d}(\theta) +  w_r \mathcal{L}_{r}(\theta) + w_b \mathcal{L}_{b}(\theta),
\end{equation}
where
\begin{equation}
\mathcal{L}_{r}(\theta) = \frac{1}{N_r} \sum_{n=1}^{N_r} \left| \mathcal{F}(y_{NN,r}^n,x_r^n, t_r^n;\theta) \right|^2\,,
\end{equation}
\begin{equation}
\mathcal{L}_{b}(\theta) = \frac{1}{N_b} \sum_{k=1}^{N_b} \left| \mathcal{B}(y_{NN,b}^k,x_b^k, t_b^k;\theta) \right|^2\,.
\end{equation}
Here, $(x_r^n, t_r^n)$, $n=1,\ldots,N_r$, and $(x_b^k, t_b^k)$, $k=1,\ldots,N_b$, are scattered points within the domain $\Omega$ and on the boundary $\partial \Omega$, respectively, called \textit{residual points} at which the PINN is asked to satisfy the known physics, while ${y}_{NN,r}^n= {y}_{NN}(x_r^n, t_r^n)$ and ${y}_{NN,b}^k = {y}_{NN}(x_b^k, t_b^k)$. Note that, this time, the subscripts ``$r$" and ``$b$" indicate the relation with the residual of the PDE and the boundary condition, respectively.
The additional residual terms, $\mathcal{L}_{r}$ and $\mathcal{L}_{b}$, quantify the discrepancy of the neural network surrogate $y_{NN}$ with the underlying differential operator and its initial or boundary conditions, respectively. Finally, $w_{r}$, $w_{b}$, and $w_{d}$ are the weights associated with each contribution, which might be higher than one. In this additional process, gradients of the DNN’s output to its input (which are needed to evaluate the operators $\mathcal{F}$ and $\mathcal{B}$) are computed at each residual point recurring to automatic differentiation \cite{baydin2017}.
The PINN so constructed can predict the solution even far from the experimental data points, behaving much better than a standard neural network \cite{Bertaglia2022b}.

At this point, the problem of training the neural network corresponds to solving a stochastic optimization problem, which is typically addressed using a stochastic gradient descent (SGD) algorithm, such as the Adam advanced optimizer \cite{Adam}.
After finding the optimal set of parameter values $\theta^*$ by minimizing the loss function \eqref{eq:general-nn-loss}, i.e.,
\begin{equation}
\theta^* = \mathrm{argmin}\, \mathcal{L}(\theta),
\end{equation}
the neural network surrogate $y_{NN}(x,t; \theta^*)$ can be evaluated at any point in the domain to get the solution to the desired problem.

Note that in the presented framework we have tacitly assumed that the $\mathcal{F}$ operator is completely known. However, often the physical equations governing the dynamics under consideration are only partially known. For example, the model parameters may not be easily assessed. We would then be faced with solving not only a forward problem of reconstruction and prediction, but also an inverse problem of finding unknown coefficients.
In the context of inverse problems, the unknown set of $n_p$ physical (scalar) parameters $\xi = (\xi_1,\ldots,\xi_{n_p}) \in \mathbb{R}^{n_p}$ is treated as a set of learnable parameters exactly like $\theta$, and, similarly for weights and biases, we search for those coefficients $\xi$ that best fit the available information. As a result, the training process involves optimizing $\theta$ and $\xi$ jointly:
\begin{equation}\label{eq:optim-loss}
    (\theta^*, \xi^*) = \mathrm{argmin}\,  \mathcal{L}(\theta,\xi).
\end{equation}

\section{PINN extension for optimal control problems}
\label{sec:PINN-OCP}
In this section, we define the control problem we want to tackle and then explain our method based on the PINN approach to solve it. We underline that the Optimal Control Problem -- Physics Informed Neural Network (OCP-PINN) architecture proposed here has been developed to solve simultaneously the problem of discovering partially unknown PDEs and the control problem, namely Eq. \eqref{eq:optim-loss}.

\subsection{PDE-constrained optimization problem}
Let us consider that the differential operator $\mathcal{F}$ defining the state equation presented in the previous section is now controlled by the function $u(x,t)$. We can then re-write the operator as 
\begin{equation}
\label{eq:state}
\mathcal{F}_c(y,u,x,t; \xi):=\frac{\partial y}{\partial t} - N(y,u,x,t;\xi) = 0, \quad  (x, t) \in \Omega,
\end{equation}
where we recall that $\xi$ is the set of model parameters that might be unknown a-priori, $N$ represents a generic non-linear differential operator in space, and we indicate the controlled state variable with $y$. The state equation \eqref{eq:state} is usually equipped with boundary conditions and initial conditions to guarantee existence and uniqueness, addressed in $\mathcal{B}_c(y,u;\xi)=0$ with $(x,t) \in \partial \Omega$. The solution to equation \eqref{eq:state} belongs to $L^2([0,T];V)$, where $V=H^k(Q)$ with $k$ depending on the maximum order of the derivatives of $y$ appearing in operator $N$. We refer to e.g. \cite{T10} for more details on control of PDEs.

We aim at finding the optimal control \( u(x,t) \) that, in this work, minimizes a quadratic cost functional with terminal cost, given by:
\begin{equation}
\label{cost}
J(y, u) = \frac{1}{2}\int_{0}^{T} \left(  \| y - \bar{y} \|_{L^2(Q)}^2 + \alpha\| u \|^2_{L^2(Q)} \right) dt + \frac{\alpha_T}{2} \| y_T - y_{f} \|^2_{L^2(Q)}\,,
\end{equation}
where \( \bar y \) is the desired state trajectory, $y_{f}$ is the desired final state, at the final time \( T \), \( y_T=y(x,T) \), and \( \alpha,\alpha_T>0 \) are two weighting parameters.

The resulting PDE-constrained optimization problem reads:
\begin{equation}\label{ocp}
   \mbox{find } \min_{u \in \mathcal{U}} J(y,u) \mbox{ such that } y \mbox{ satisfies } \eqref{eq:state} \,,
\end{equation}
where $\mathcal{U}$ is the space of admissible controls. In our setting, following \cite{T10} and to guarantee square-integrability we assume that $\mathcal{U}=L^2([0,T];L^2(Q))$.
The control problem \eqref{ocp} can be solved in an open-loop fashion using Lagrangian multipliers. Therefore, we define the Lagrangian function
\begin{equation}
\mathfrak{L}(y,u,p,x,t;\xi) = J(y,u) + \int_\Omega p(x,t) \cdot \mathcal{F}_c(y,u,x,t; \xi) \,dx \, dt\,,
\end{equation}
where $p(x,t)$ is the adjoint variable of the method. Indeed, from the above Lagrangian, we retrieve the so-called \textit{adjoint equation} given by the following backward linear equation:
\begin{align}\label{eq:adjoint}
\begin{aligned}
\mathcal{A}(y,p,x,t;\xi):= (y-\bar{y}) -\frac{\partial p}{\partial t} -p\,\frac{\partial N}{\partial y} &= 0, \quad  &(x, t) \in \Omega\,, \\ 
    \mathcal{B}_{\mathcal{A}}(p, x, t; \xi) &= 0, \quad &(x, t) \in \partial \Omega\,,
\end{aligned}
\end{align}
where boundary and terminal conditions for the adjoint variable, $p(x,t)$ for $x\in\partial \Omega$ and $p(x,T)$, respectively, are shortly included within the term $\mathcal{B}_\mathcal{A}$. Typically, the adjoint variable belongs to the same functional space of the state variable (see \cite{T10}).
Finally, the control \( u(x,t) \) has to satisfy the maximum principle (optimality condition), which states:
\begin{equation}
\label{eq:opt}
\mathcal{O}(p,u,x,t):=\alpha\, u  +p\,\frac{\partial N}{\partial u}=0, \quad  (x, t) \in \Omega\,.
\end{equation}

The solution of the OCP \eqref{ocp}, following the Lagrange multipliers approach, can be obtained by solving system \eqref{eq:state}-\eqref{eq:adjoint}-\eqref{eq:opt}.
For more details on the method, we invite the reader to refer, for instance, to \cite{HPUU09,T10}.

\subsection{OCP-PINNs}\label{sec:a-pinn}
Here, we propose to solve the optimal control problem and also to correctly identify possible unknown parameters of the investigated model via the same online stage using ad-hoc constructed PINNs, which we shall call OCP-PINNs.

\begin{figure}[tp!]
    \centering    
    \includegraphics[width=\textwidth]{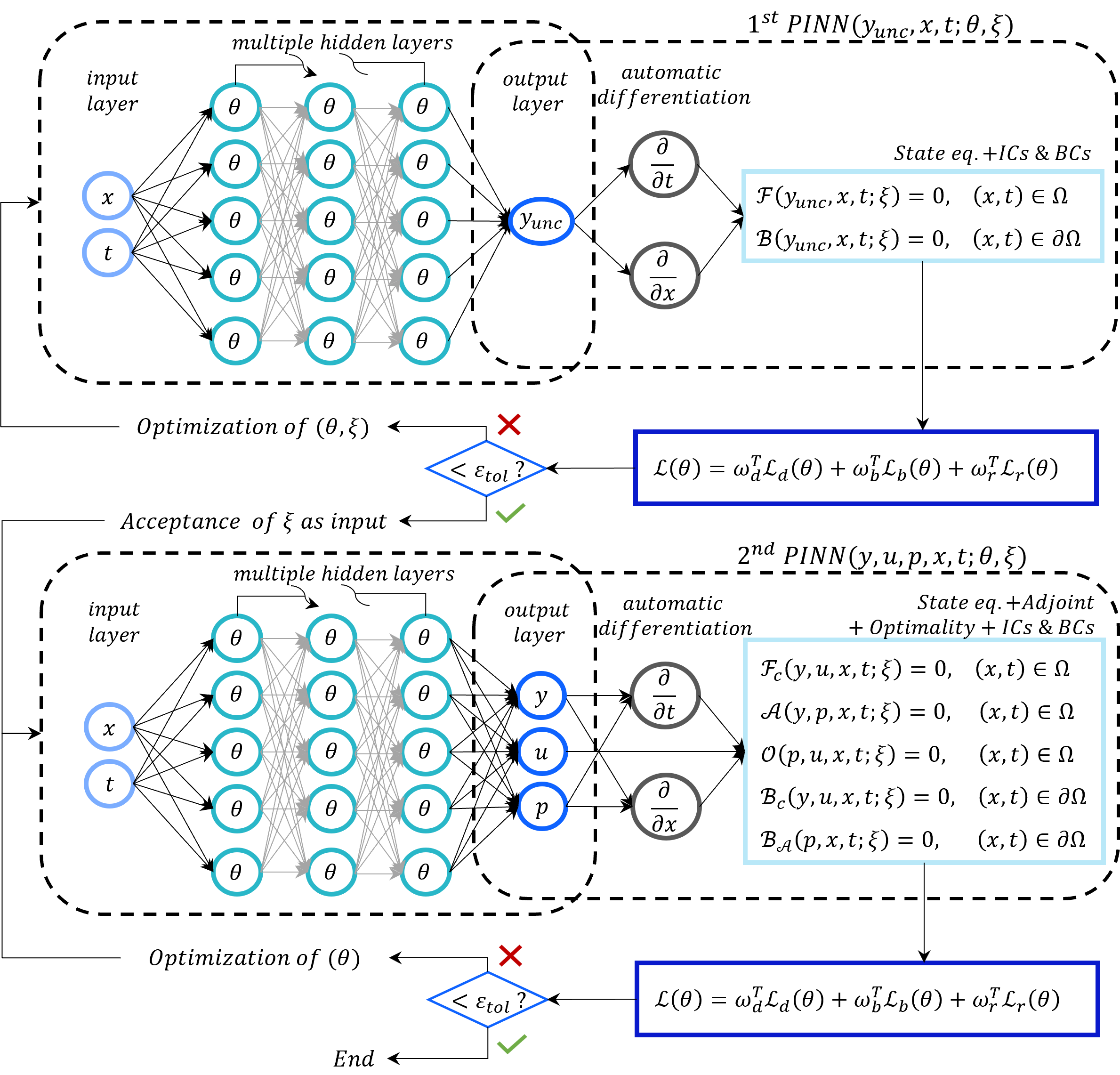}    
    \caption{OCP-PINN schematic workflow. \first{We consider an architecture composed by 2 PINNs in series: the first PINN only deals with the resolution of the inverse problem of estimating possible unknown parameters through the learning of the uncontrolled solution, while the second PINN accounts for the control problem, taking as input also the parameters discovered from the first neural network. In the latter,} the PINN architecture is integrated with the physical knowledge of the adjoint equation and the optimality condition, through the inclusion of the differential operators obtained following the Lagrange multipliers approach, as well as boundary conditions of the corresponding variables.}
    \label{OCP-PINN-scheme}
\end{figure}

Starting from the basic framework presented in Section \ref{sec:NN-PINN}, \first{we consider an architecture composed by 2 PINNs in series: the first PINN only deals with the resolution of the inverse problem of estimating possible unknown parameters through the learning of the uncontrolled solution (hence, the output of this neural network is only the uncontrolled state, $y_{unc}$), while the second PINN accounts for the control problem, taking as input also the parameters discovered from the first PINN, having as outputs the control variable $u$ and the adjoint state $p$, together with the controlled state $y$. Therefore, the loss function of the first PINN reads exactly as in \eqref{eq:general-pinn-loss}, being 
\[\mathcal{L}_{r}(\theta,\xi) = \frac{1}{N_r} \sum_{n=1}^{N_r} \left| \mathcal{F}(y_{unc,NN,r}^n, x_r^n, t_r^n;\theta,\xi) \right|^2\,.\]
While in the second PINN,} the loss function term accounting for the physical knowledge of the dynamics, $\mathcal{L}_r$ in \eqref{eq:general-pinn-loss}, is extended to include the additional information given by the setting of the control problem obtained using Lagrangian multipliers. Particularly, we recall the loss term accounting for the residual of the state equation,
\begin{equation}
\mathcal{L}_{r,\mathcal{F}}(\theta) = \frac{1}{N_r} \sum_{n=1}^{N_r} \left| \mathcal{F}_c(y_{NN,r}^n, u_{NN,r}^n, x_r^n, t_r^n;\theta) \right|^2\,,
\end{equation}
and add a loss term to account for the residual for the adjoint equation and another for the optimality condition:
\begin{equation}
\mathcal{L}_{r,\mathcal{A}}(\theta) = \frac{1}{N_r} \sum_{n=1}^{N_r} \left| \mathcal{A}(y_{NN,r}^n, p_{NN,r}^n, x_r^n, t_r^n;\theta) \right|^2\,,
\end{equation}
\begin{equation}
\mathcal{L}_{r,\mathcal{O}}(\theta) = \frac{1}{N_r} \sum_{n=1}^{N_r} \left| \mathcal{O}(p_{NN,r}^n, u_{NN,r}^n, x_r^n, t_r^n;\theta) \right|^2\,,
\end{equation}
so that 
\begin{equation}\label{eq:residual-adjoint}
\mathcal{L}_r = \mathcal{L}_{r,\mathcal{F}} + \mathcal{L}_{r,\mathcal{A}} + \mathcal{L}_{r,\mathcal{O}}\,.
\end{equation}
Similarly, the loss term related to boundary conditions (recall that here by boundary conditions we refer to both the usual boundary conditions and the initial and final boundary conditions since we treat the space-time domain as a single 2D space) will read 
\begin{equation}\label{eq:residualBC-adjoint}
\mathcal{L}_{b} = \mathcal{L}_{b,\mathcal{F}} + \mathcal{L}_{b,\mathcal{A}},
\end{equation}
where
\begin{equation}
\mathcal{L}_{b,\mathcal{F}}(\theta;\xi) = \frac{1}{N_b} \sum_{k=1}^{N_b} \left| \mathcal{B}_{c}(y_{NN,b}^k,u_{NN,b}^k,x_b^k, t_b^k;\theta,\xi) \right|^2 \,,
\end{equation}
\begin{equation}
\mathcal{L}_{b,\mathcal{A}}(\theta;\xi) = \frac{1}{N_b}\sum_{k=1}^{N_b} \left| \mathcal{B}_\mathcal{A}(p_{NN,b}^k,x_b^k, t_b^k;\theta,\xi) \right|^2\,.
\end{equation}
Notice that $\mathcal{L}_{b,\mathcal{F}}$ accounts for the boundary conditions of both the state variable and the control. For a schematic representation of the discussed OCP-PINN architecture, the reader can refer to Fig. \ref{OCP-PINN-scheme}.

In the above setting, all possible residual terms that may be known about the physics of the observed phenomenon are considered. However, it is not certain that all of them can actually be defined a priori, and it is in this very challenging setting that we will focus on the numerical tests that follow in the next section. 
In this paper, we are interested in observing the behavior of the OCP-PINNs structured here in case studies where the initial and boundary conditions of the problem are only available at a sparse level in the computational domain and the governing physical laws are only partially known, solving Eq. \eqref{eq:optim-loss}. \first{Specifically, in the numerical tests to follow we will consider only a few scattered data available in the space-time domain $\Omega$ of the uncontrolled solution and very few data only of the initial condition of the controlled state variable. It is emphasized that the information about the initial condition of the controlled problem actually could be recovered directly from the reconstruction obtained from the first PINN, related to the uncontrolled problem, since the initial state is clearly the same. This would mean, in essence, to assume that no data on the controlled problem is available at all.}

\begin{remark}
We remark here that, following the above-presented procedure and appropriately embedding the physical residual (as discussed in \cite{Bertaglia2022b,Bertaglia2022a}), it is straightforward to construct Asymptotic-Preserving Neural Networks (APNNs) to solve optimal control problems involving multiscale dynamics, namely OCP-APNNs for which the loss function of the full model-constraint converges to the loss function of the corresponding asymptotic reduced-order model. See also \cite{Jin2022} for further details on APNNs. We leave the treatment of these kinds of problems to future investigations.
\end{remark}

\begin{remark}
\first{It is worth to mention that the question of whether to use a single neural network with different outputs (in our case, the 3 outputs $y$, $u$, and $p$), or alternatively, to train separate neural networks for each physical variable, remains debated. The fundamental distinction lies in whether or not the network representations of the various quantities should incorporate shared weights. The proposed configuration represents a balanced approach that has proven to be relatively easier to train than others for the numerical tests considered in this work.}
\end{remark}

\section{Numerical tests}
\label{sec:tests} 
In this section, we test the performance of the proposed OCP-PINN approach through various numerical experiments in conditions in which available data are very sparse and the physical knowledge of the dynamics is incomplete. For each test case, the state and adjoint equations as well as the optimality condition are presented. We specify here that in the following we will always consider a cost functional \eqref{cost} in which  $\bar{y} = y_f = 0$; hence we will always aim to direct the dynamics toward the zero solution. 

In all the tests that follow, for the sake of simplicity, we will indicate the reference solution of the state variable with $y(x,t;u^*)=y(u^*)$, where $u^*$ is the reference optimal control. This solution is calculated using the standard gradient method \cite{HPUU09} and will provide the synthetic dataset from which to select some points to feed the neural network in the training phase. The reference solution is then used for comparison with the results obtained by the proposed OCP-PINN in terms of the state variable output to the network, denoted by $y_{PINN}$. Furthermore, to perform error analysis, we indicate with $u_{PINN}$ the control variable obtained from the proposed OCP-PINN approach as an output of the neural network, while we define $y(x,t;u_{PINN})=y(u_{PINN})$ the numerical approximation of the state equation resulting from plugging the PINN control $u_{PINN}$ directly into the PDE and solving numerically the evolutionary problem for the fixed control. These quantities will permit us to further evaluate the quality of our results. In particular, we define the following $L^2$ relative norms that will be used in the following:
$$\mathcal{E}_1:=\frac{\|y(u^*)-y_{PINN}\|_{L^2([0,T];L^2(Q))}}{\|y(u^*)\|_{L^2(0,T];L^2(Q))}},\,\, \mathcal{E}_2:=\frac{\|y(u^*)-y(u_{PINN})\|_{L^2([0,T];L^2(Q))}}{\|y(u^*)\|_{L^2([0,T];L^2(Q))}},$$
$$\mathcal{E}_3:=\frac{\|u^*-u_{PINN}\|_{L^2([0,T];L^2(Q))}}{\|u^*\|_{L^2([0,T];L^2(Q))}}, \quad \mathcal{E}_4:=\frac{\|y(u_{PINN})-y_{PINN}\|_{L^2([0,T];L^2(Q))}}{\|y(u_{PINN})\|_{L^2([0,T];L^2(Q))}}.$$

We also report the values of the cost functional computed using the reference solution, $J(y(u^*),u^*)$, the value of the cost resulting from the use of PINN,  $J(y(u_{PINN}),u_{PINN})$, and the relative error between them
\[\mathcal{E}_5 := \frac{\left|J(y(u^*),u^*)-J(y(u_{PINN}),u_{PINN})\right|}{J(y(u^*),u^*)}.\]
	
	 In all simulations, a 5-layer feed-forward neural network is used, having 3 hidden layers with 64 nodes each. We recur to the Adam method \cite{Adam} for the optimization process, fixing a learning rate $\mathrm{LR}=10^{-3}$ and considering the \texttt{tanh} as activation function. Finally, in the loss function \eqref{eq:general-pinn-loss} we set the following weights: $w_d=w_b=10$ and $w_r = 1$.

\subsection{Test 1: Viscous Burgers' equation}\label{sec:test1}
The first test concerns the control of the viscous Burgers' equation, namely:
\begin{equation}
    \frac{\partial y}{\partial t} + y\frac{\partial y}{\partial x} = \nu \frac{\partial^2 y}{\partial x^2}  + u,
\end{equation}
in which we recall that $y=y(x,t)$ is the state variable and $u=u(x,t)$ is the control variable. Here $\nu>0$ identifies the viscosity coefficient. We consider a spatio-temporal domain $\Omega = [a,b]\times[0,T]$, with $a=-4$, $b=4$, $T=4$, and the following initial and boundary conditions, respectively:
\begin{align}\label{burg}
\begin{aligned}
y(x,0) &= \frac{e^{-0.5x^2}}{\sqrt{2\pi}},\qquad &&x\in[a,b],\\
y(a,t) &= y(b,t) = 0,\qquad &&t\in[0,T] .
\end{aligned}
\end{align}
For the control problem, in the cost functional \eqref{cost}, we fix $\bar{y}=y_f=0$, $\alpha=10$, and $\alpha_T=0.1$, hence we have
\begin{equation}
J(y, u) = \frac{1}{2}\left(\int_{0}^{T} \left(  \| y \|_{L^2([a,b])}^2 + 10\| u \|^2_{L^2([a,b])} \right) dt +  0.1\| \,y_T\|^2_{L^2([a,b])}\right)\,.
\end{equation}
Following the Lagrangian multipliers approach discussed in Section \ref{sec:a-pinn}, the adjoint equation for the variable $p=p(x,t)$ reads
\begin{equation}
    y -\frac{\partial p}{\partial t} -y\frac{\partial p}{\partial x} = \nu \frac{\partial^2 p}{\partial x^2}.
\end{equation}
This equation results equipped with the following boundary and terminal conditions:
\begin{align}
\begin{aligned}
p(a,t) &= p(b,t) = 0,  \qquad &&t\in[0,T],\\
p(x,T) &= - 0.1\,y(x,T),  \qquad &&x\in[a, b].
\end{aligned}
\end{align}
Finally, the optimality condition results
\begin{equation}\label{eq:opt-burgers}
  10\,u = p.  
\end{equation}
Note that the above optimality condition is rather simple in this specific case, and we immediately observe that the control $u$ gets the same boundary conditions as the adjoint $p$. 

In this numerical test, we aim at solving the control problem while also discovering the unknown viscosity coefficient $\nu$ in \eqref{burg} using the proposed OCP-PINN approach. In particular, we consider two scenarios, one in which the parameter to be discovered is $\nu = 0.5$ (case a), i.e., in which the problem is very diffusive, and the second in which $\nu = 0.05$ (case b), i.e., a case of milder diffusion and dominance of the advective behavior of the solution.

In both cases, we train the OCP-PINN using a scattered dataset of only $N_d=10$ training points given for \first{both uncontrolled and controlled state} variables, selected from the reference solution, \first{but considering only initial condition data for the controlled one. The knowledge of the uncontrolled data will allow the inverse problem to be solved}. No \first{additional} data is considered available, making the problem very challenging for the neural network. Concerning boundary conditions, we impose them at a sparse level, in $N_b=18$ points \first{for all variables}, and we consider available terminal conditions \first{for the adjoint $p$}. Finally, we fix $N_r=336$ residual nodes uniformly distributed in $\Omega$ for evaluating the physical loss terms.
The training is executed for a number of iterations (epochs) to satisfy a tolerance on the error given by the loss function \eqref{eq:general-pinn-loss}-\eqref{eq:residual-adjoint}-\eqref{eq:residualBC-adjoint} of less than \first{$\varepsilon_{tol}=10^{-6}$}.

Notice that since the optimality condition \eqref{eq:opt-burgers} is given by an algebraic equation, for this test we could reduce the number of outputs of the \first{second-in-series PINN} to 2 and directly evaluate $p$ as a function of $u$ through \eqref{eq:opt-burgers} (or vice-versa). However, for this test, we chose to keep the generic architecture of the proposed OCP-PINN, with the 3 outputs as presented in Fig. \ref{OCP-PINN-scheme}.

\begin{figure}[tp!]
    \centering    
    \includegraphics[width=0.7\textwidth]{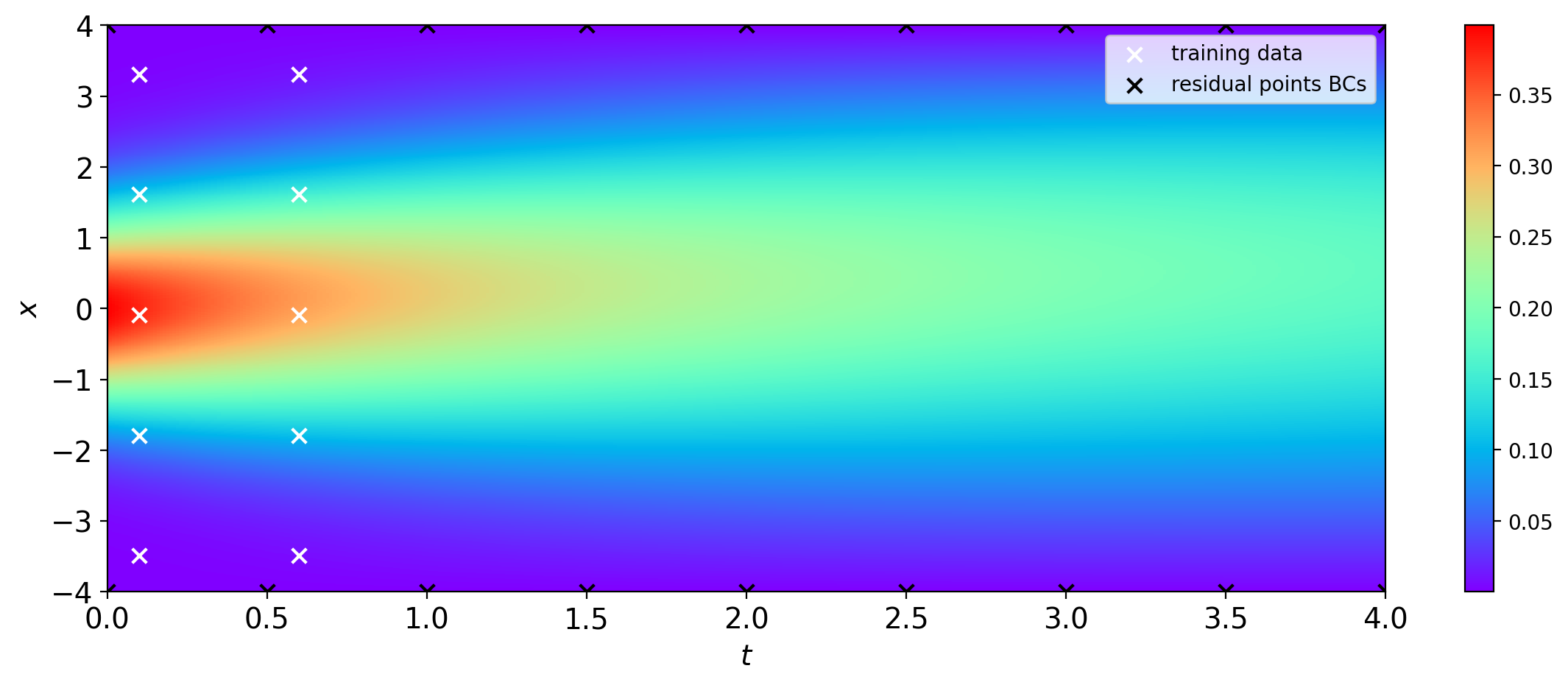}    
    \includegraphics[width=0.7\textwidth]{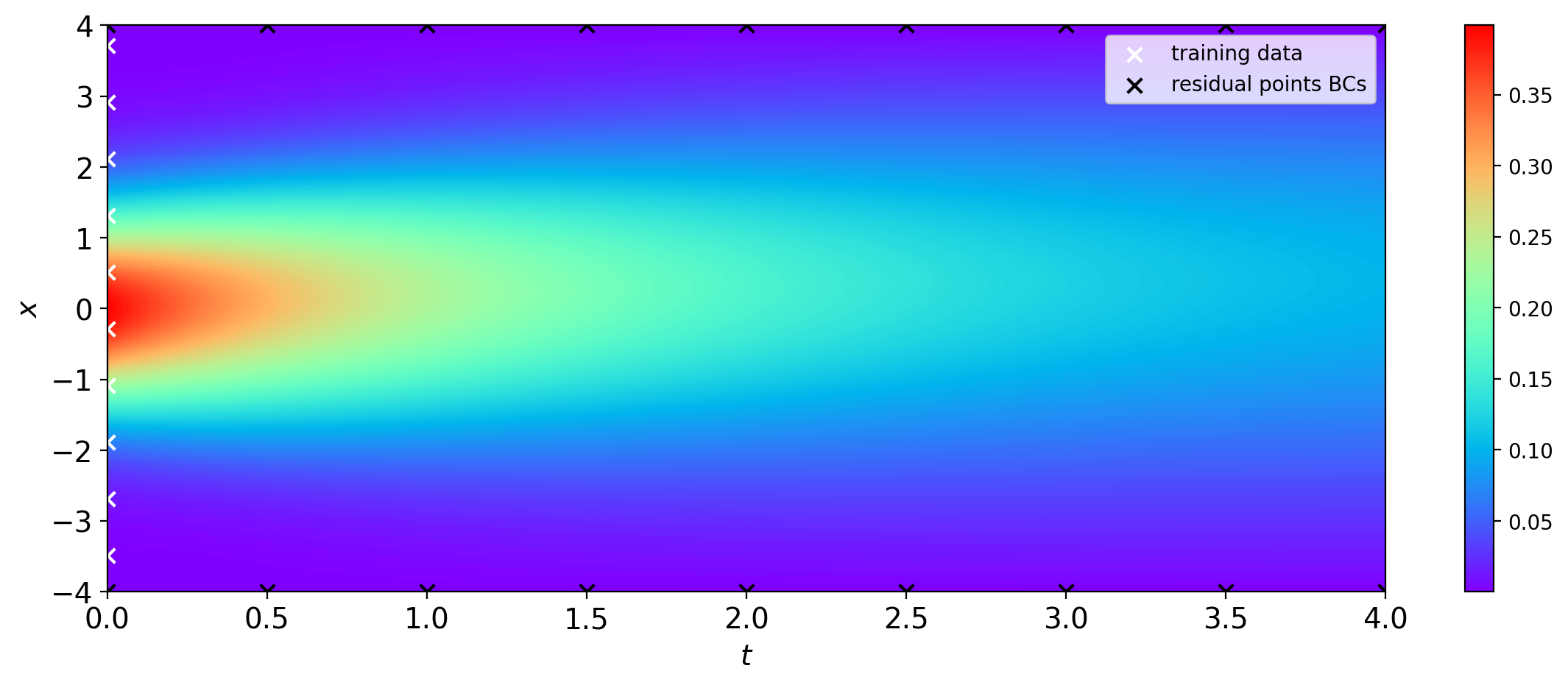}       
    \caption{Test 1(a): Burgers' equation with $\nu=0.5$. Data points of the \first{uncontrolled (top) and controlled (bottom)} state variables used for the training of the OCP-PINN.}
    \label{coll-11}
\end{figure}

\begin{figure}[tp!]
    \centering
    \includegraphics[width=0.48\textwidth]{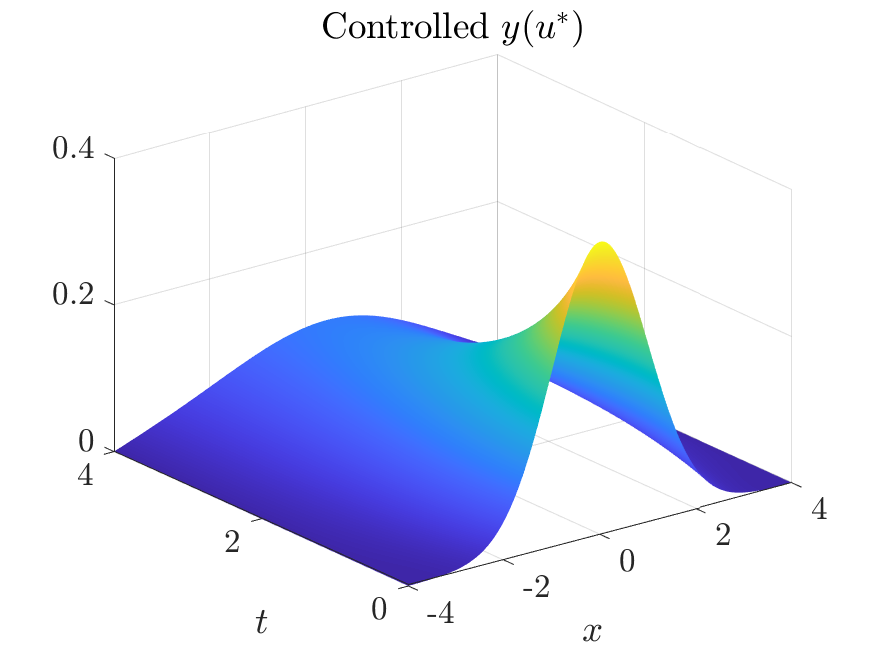}
    \includegraphics[width=0.48\textwidth]{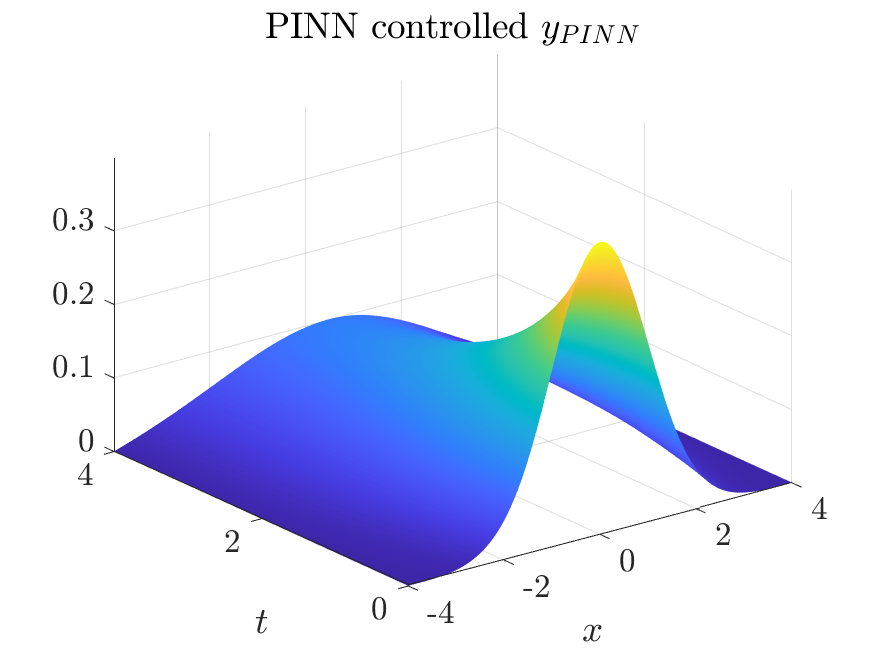}
    \includegraphics[width=0.48\textwidth]{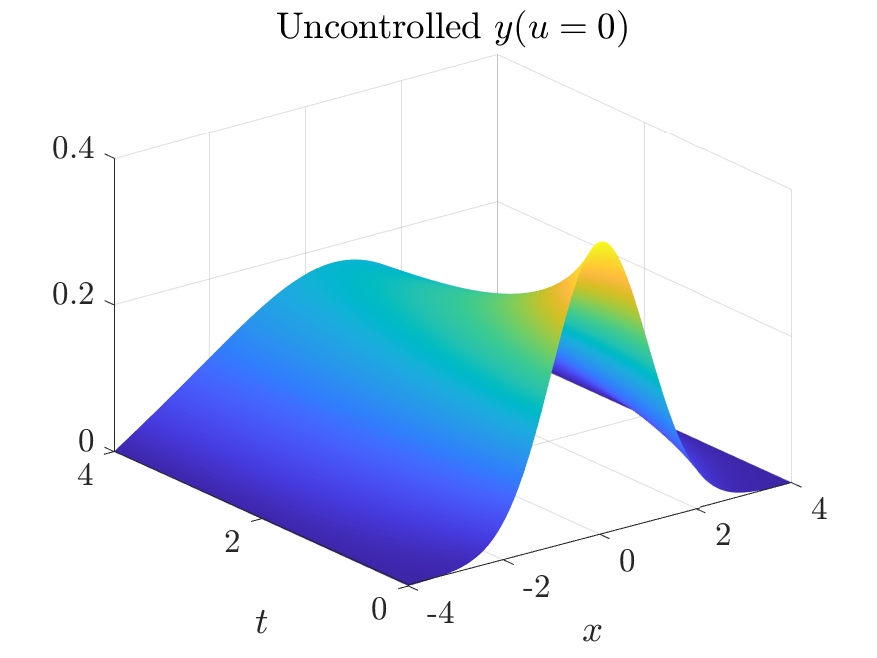}
    \caption{Test 1(a): Burgers' equation with $\nu=0.5$. Reference controlled solution of the state variable (top left), OCP-PINN controlled solution (top right), and reference uncontrolled solution (bottom).}
    \label{fig:burg1}
\end{figure}

\begin{figure}[tp!]
    \centering
    \includegraphics[width=0.48\textwidth]{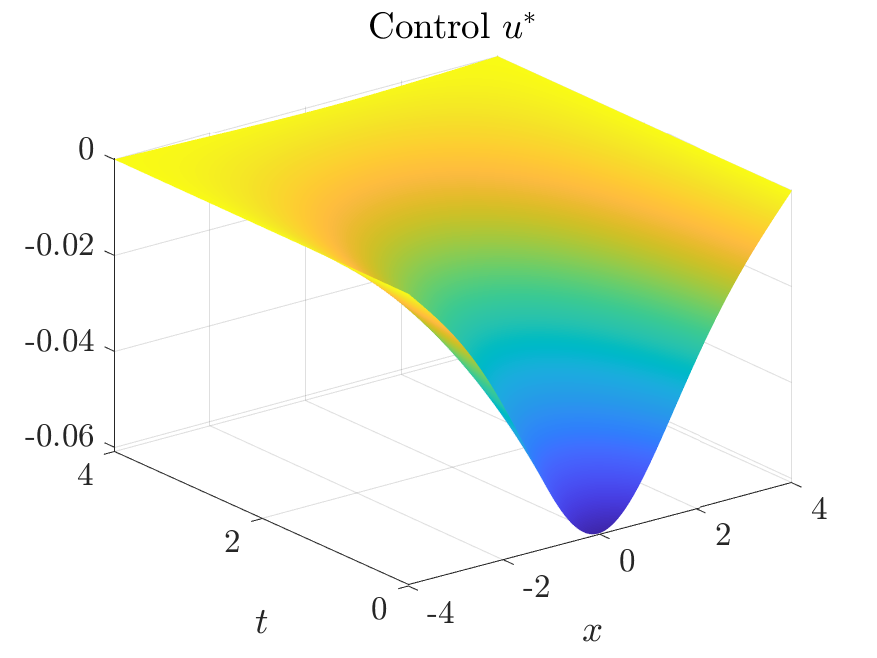}
    \includegraphics[width=0.48\textwidth]{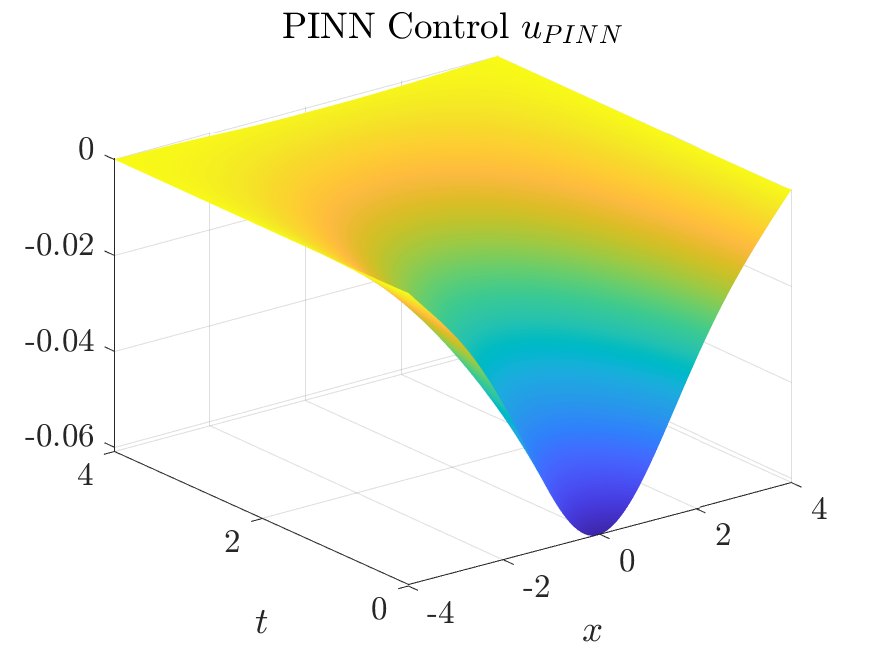}
    \caption{Test 1(a): Burgers' equation with $\nu=0.5$. Reference control variable (left) and control obtained by applying the OCP-PINN (right).}
    \label{burg:contr}
\end{figure}

\begin{figure}[tp!]
    \centering
    \includegraphics[width=0.48\textwidth]{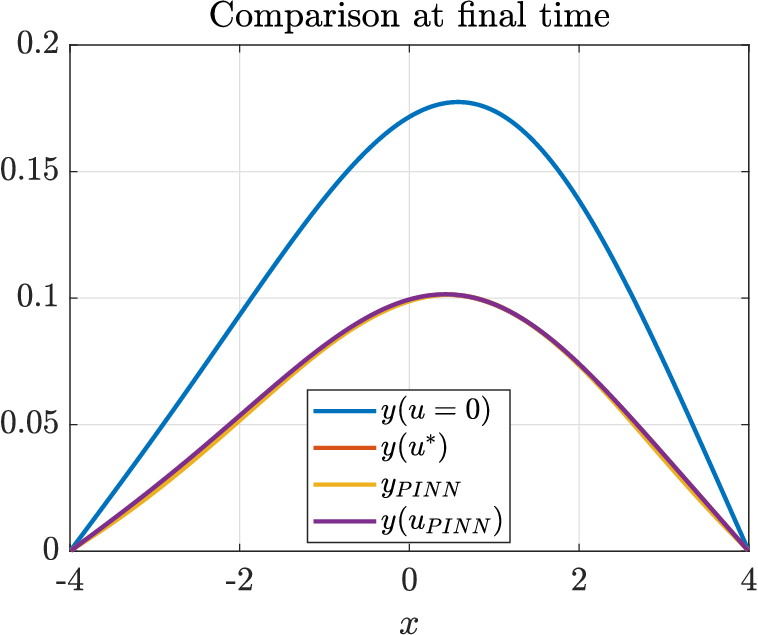}    \hspace{0.25cm}
    \includegraphics[width=0.48\textwidth]{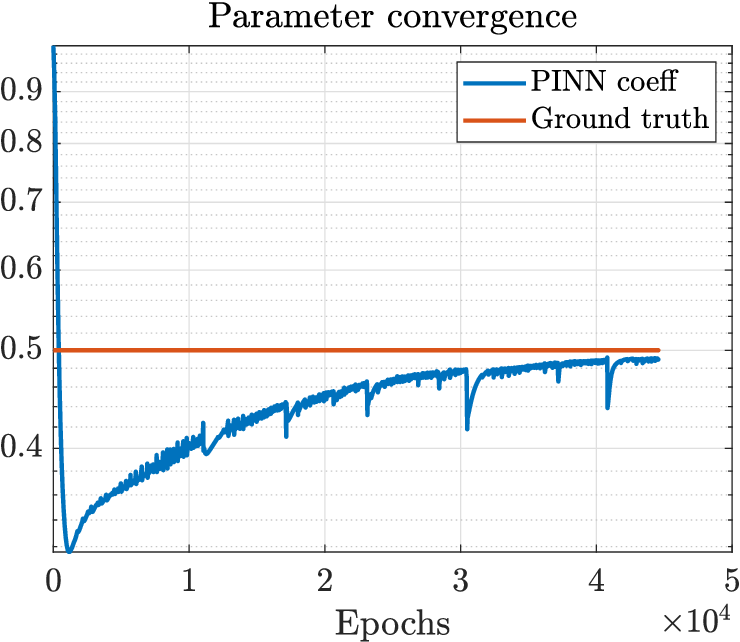}
    \caption{Test 1(a): Burgers' equation with $\nu=0.5$. Comparison of the solutions at final time $T=4$ (left) and convergence history for discovering the unknown parameter $\nu$ (right).}
    \label{burg:comp}
\end{figure}

\paragraph{Test 1(a): the case $\nu=0.5$.}

We show the training points $N_d$ and the residual points $N_b$ used for the training phase in Fig. \ref{coll-11}. 
In Fig. \ref{fig:burg1}, we show the controlled solution obtained with the OCP-PINN, $y_{PINN}$, compared with the reference controlled, $y(u^*)$, and uncontrolled solutions. It is visually clear how the controlled solution tends to reach the zero state faster than the uncontrolled one, which is the goal of the OCP investigated. We can affirm that these results confirm the validity of the approach. Similar considerations can be made by looking at Fig. \ref{burg:contr}, where we present the reference solution for the control variable $u^*$ in comparison with the PINN control $u_{PINN}$. Indeed, here we observe a very good agreement between the two solutions, \first{even at the final time, in face of} the noted difficulty of the neural network in reconstructing nearly-zero states, as expected for $u$ at the final time.
Moreover, in Fig. \ref{burg:comp}, we show in the left panel a snapshot of the solutions computed with the different techniques, together with the uncontrolled state $y_{unc}(x,t)=y(u=0)$, at the final time $T=4$. We recall that $y(u_{PINN})$ identifies the numerical approximation of the state equation resulting from plugging the PINN control $u_{PINN}$ directly into the PDE and solving numerically the evolutionary problem for the fixed control. In this plot, we can observe a minor deviation between the controlled solutions with respect to the reference controlled state.

Regarding the solution of the inverse problem for identifying the unknown parameter of the model, we show the PINN history of the search for the $\nu$ coefficient in the right panel of Fig. \ref{burg:comp}. Although an unpleasant oscillatory pattern is noticeable, the proposed PINN approach finds \first{$\nu_{PINN} = 0.4905$} starting with an initial guess of $\nu_0 = 1$, which reflects a very good identification.

With respect to the reference solution $y(u^*)$, we computed the following relative $L^2$ norms: \first{$\mathcal{E}_1= 0.007618$ and $\mathcal{E}_2 = 0.002075$}.
We observe that these relative errors are of order $O(10^{-3})$ and are quite consistent with each other.
Furthermore, it is interesting to note that the relative error between the two controls and the error between $y(u_{PINN})$ and $y_{PINN}$ is coherent with those shown above, indeed \first{$\mathcal{E}_3 =0.008515$ while $\mathcal{E}_4= 0.009527$}. Concerning the cost functional, its value computed with the reference solution is $J(y(u^*),u^*)=25.738890$, while the one obtained using PINN results in $J(y(u_{PINN}),u_{PINN})=25.738134$, and the relative error is $\mathcal{E}_5 = 2.937\cdot 10^{-5}$.


\paragraph{Test 1(b): the case $\nu = 0.05$.} 

\begin{figure}[tp!]
    \centering    
    \includegraphics[width=0.7\textwidth]{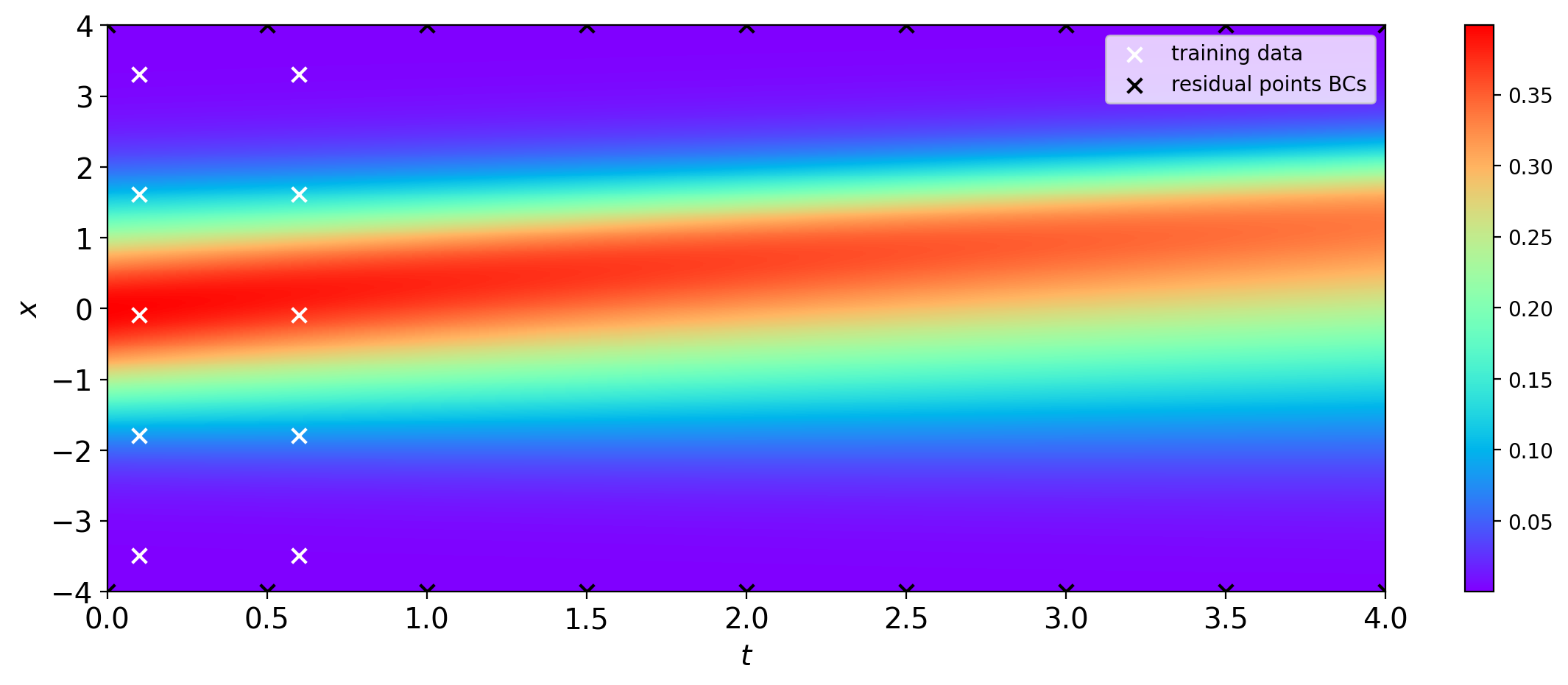}    
    \includegraphics[width=0.7\textwidth]{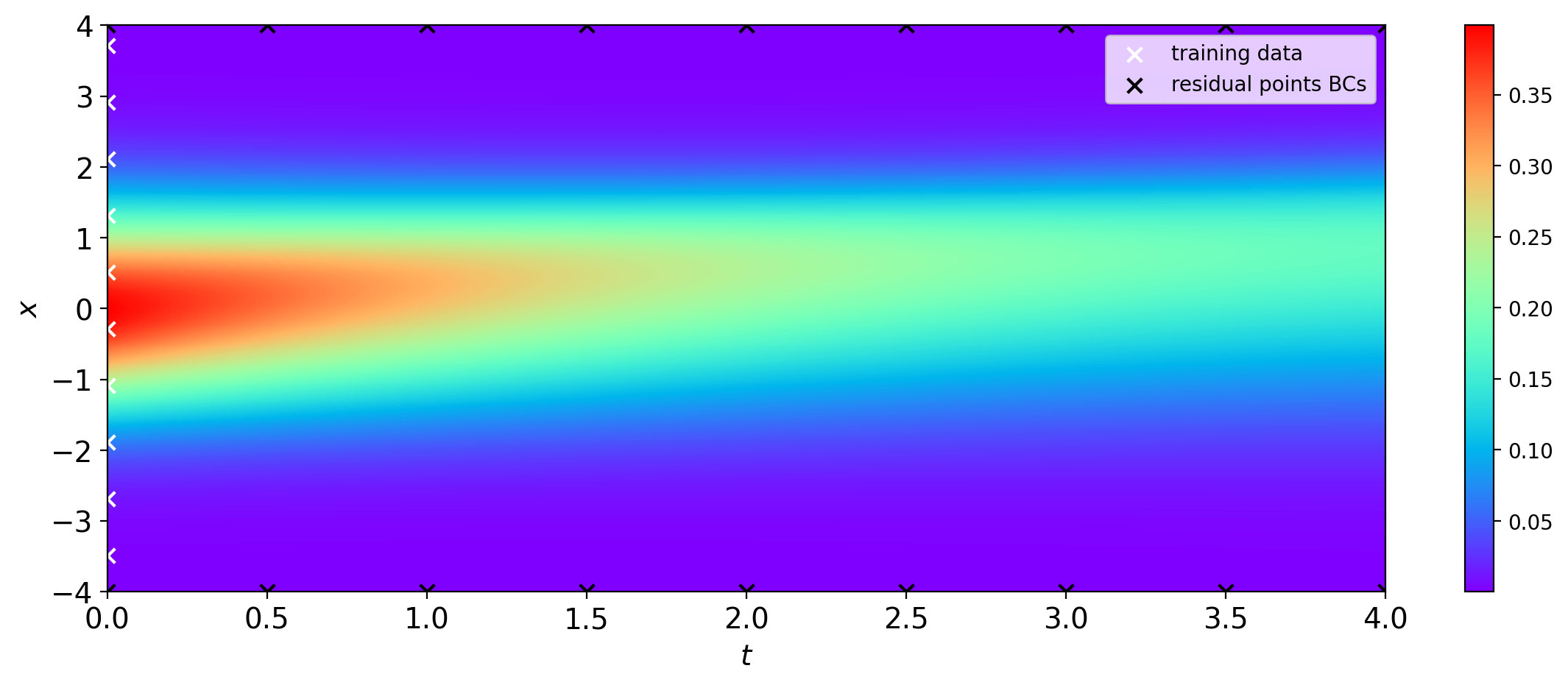}    
    \caption{Test 1(b): Burgers' equation with $\nu=0.05$. Data points of the \first{uncontrolled (top) and controlled (bottom)} state variable used for the training of the OCP-PINN.}
    \label{coll-12}
\end{figure}

\begin{figure}[tp!]
    \centering
    \includegraphics[width=0.48\textwidth]{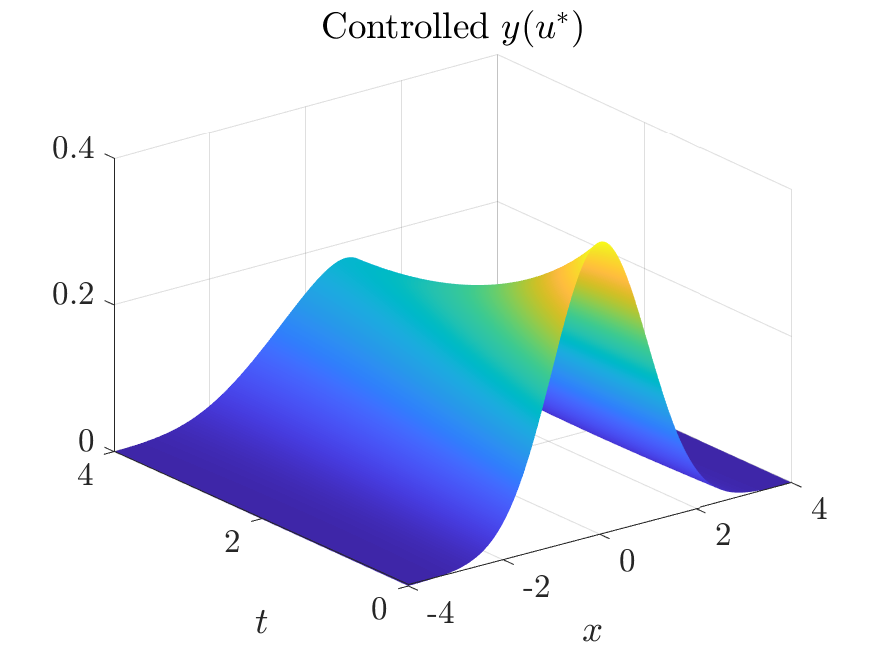}
    \includegraphics[width=0.48\textwidth]{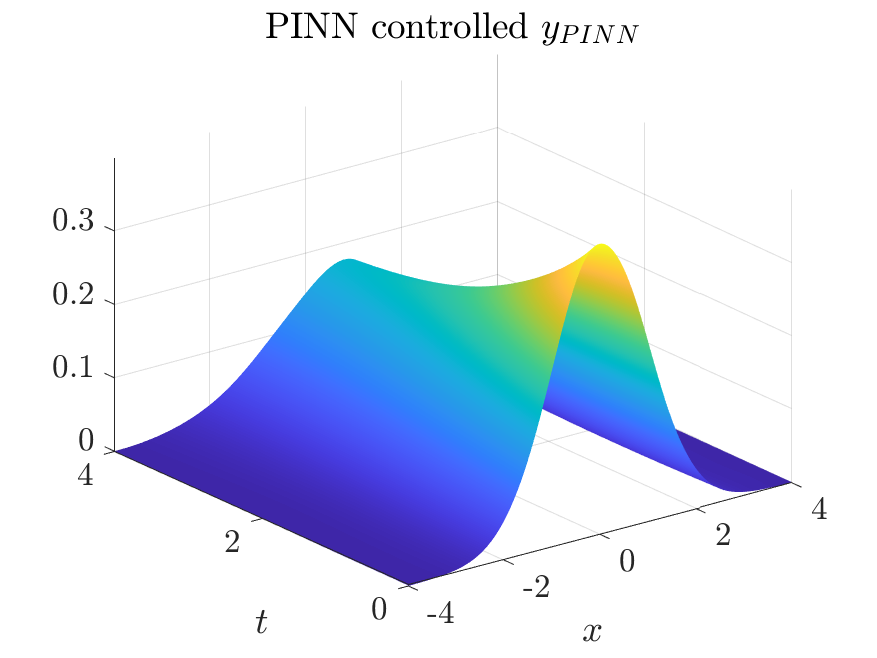}
    \includegraphics[width=0.48\textwidth]{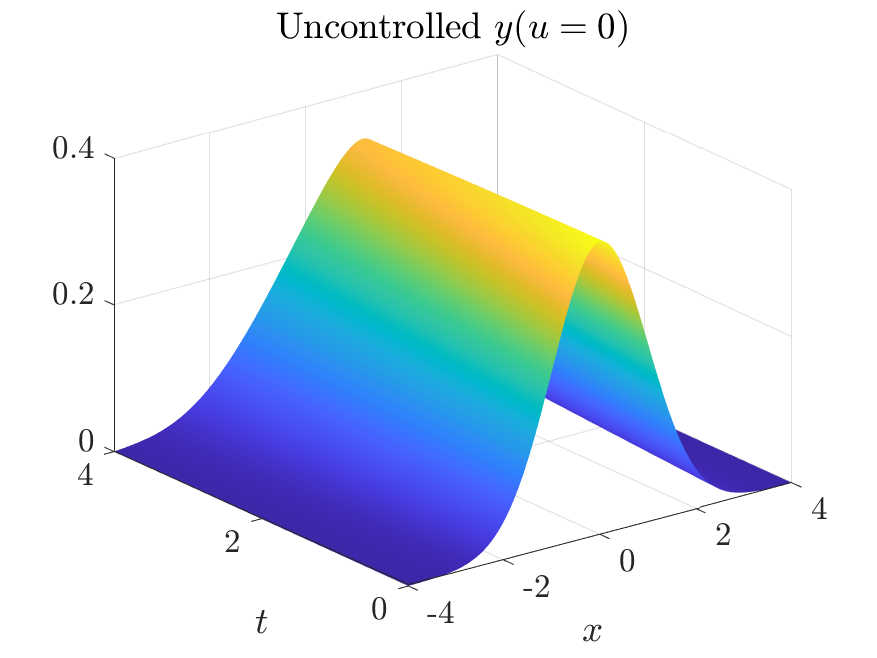}
    \caption{Test 1(b): Burgers' equation with $\nu=0.05$. Reference controlled solution of the state variable (top left), OCP-PINN controlled solution (top right), and reference uncontrolled solution (bottom).}
    \label{fig:burg1500}
\end{figure}

\begin{figure}[tp!]
    \centering
    \includegraphics[width=0.48\textwidth]{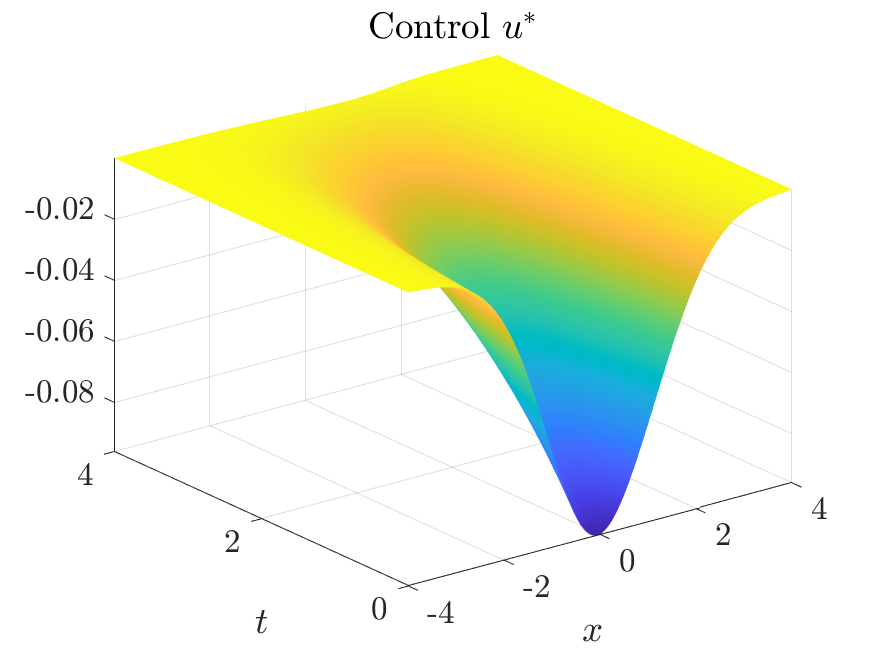}
    \includegraphics[width=0.48\textwidth]{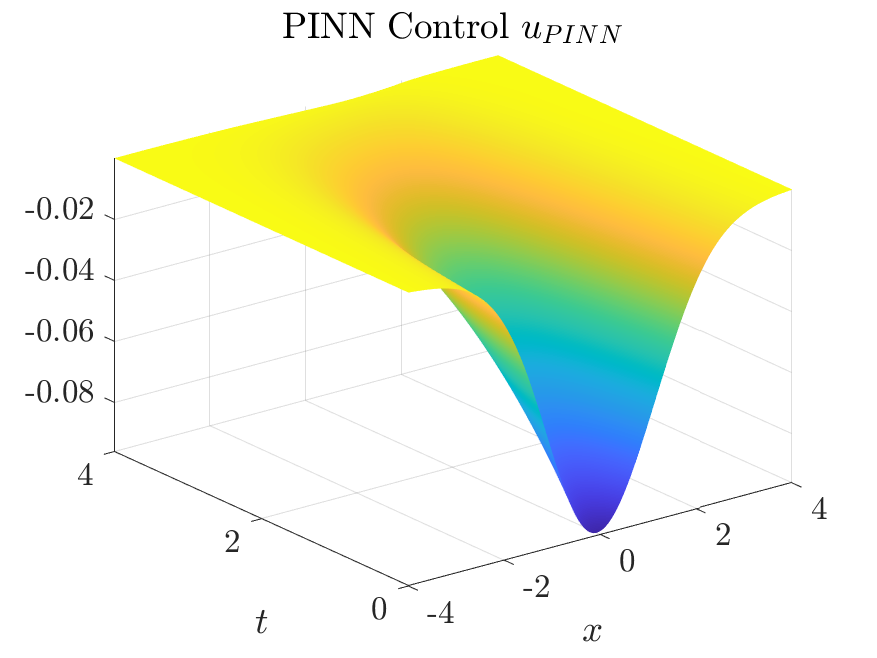}
    \caption{Test 1(b): Burgers' equation with $\nu=0.05$. Reference control variable (left) and control obtained by applying the OCP-PINN (right).}
    \label{burg:contr500}
\end{figure}

\begin{figure}[tp!]
    \centering
    \includegraphics[width=0.48\textwidth]{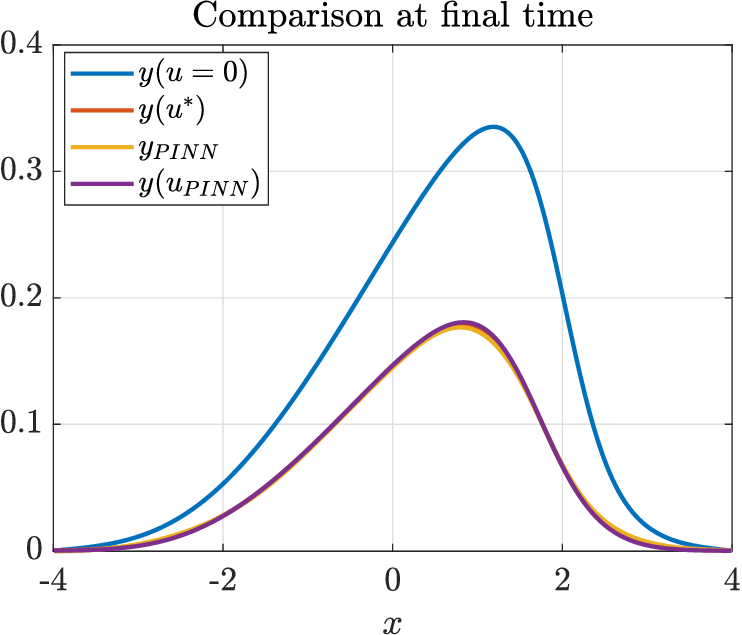}\hspace{0.25cm}
    \includegraphics[width=0.48\textwidth]{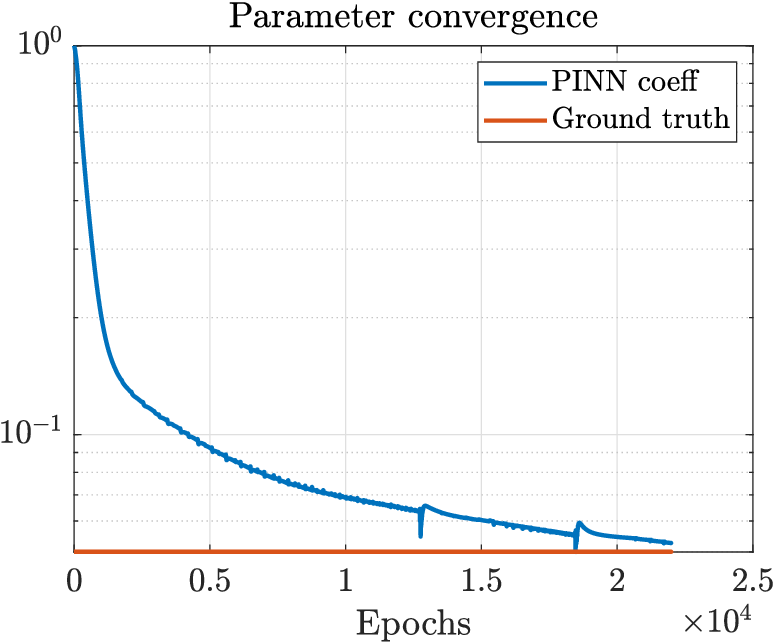}
    \caption{Test 1(b): Burgers' equation with $\nu=0.05$. Comparison of the solutions at final time $T=4$ (left) and convergence history for discovering the unknown parameter $\nu$ (right).}
    \label{burg:comp500}
\end{figure}
The location of the data points $N_d$ and $N_b$ used to train the neural network is shown in Figure \ref{coll-12}.
We present the solution of this test in terms of the state variable in Fig. \ref{fig:burg1500}. Precisely, in this figure, we compare the controlled reference state with the proposed OCP-PINN solution, which are found to be in remarkable agreement. In the same figure, we also show the uncontrolled state. In this second scenario, we are considering a smaller viscosity than in the previous case, and its effect is most evident when looking at the uncontrolled solution, which is much less diffused than in the previous example. At the same time, it appears clear how the controlled solution is going to the zero state faster than the uncontrolled one, according to the goal of the OCP. 
Figure \ref{burg:contr500} permits the reader to compare the control configuration obtained with the OCP-PINN approach with the reference one, again showing a very good agreement \first{also} concerning the shape of the solution at the final time.

In the left panel of Fig. \ref{burg:comp500}, we show the profile of different solutions at the final time $T=4$. We see how the controlled solutions $y(u^*), y(u_{PINN})$, and $y_{PINN}$ agree up to a small error on the behavior, also catching the influence of the non-linearity which makes the problem more complicated due to low viscosity. 

We show the convergence in the identification of the unknown viscosity coefficient in the right panel of the same Fig. \ref{burg:comp500}. The final value obtained is \first{$\nu_{PINN} = 0.0527$} with a very far initial guess of $\nu_0 = 1$, confirming an appropriate identification. The convergence is rather \first{fast and takes $2.1\cdot 10^4$} epochs. Furthermore, it can be seen that after \first{$10^4$} epochs we are already quite close to the desired configuration. 

The quality of the results obtained is reflected in the following relative $L^2$ norms:
\first{$\mathcal{E}_1= 0.007571$, $\mathcal{E}_2= 0.004782 $, $\mathcal{E}_3 =0.010782$, and $\mathcal{E}_4= 0.011571 $}, which are of order $O(10^{-2})$-$O(10^{-3})$. The value of the cost functional computed with the reference solution is $J(y(u^*),u^*)=35.641479$, while for the cost functional obtained using PINN we have $J(y(u_{PINN}),u_{PINN})=35.639545$, with relative error $\mathcal{E}_5 = 5.426262\cdot 10^{-5}$.


\subsection{Test 2: Allen-Cahn equation}\label{sec:test2}

The second test concerns the control of the Allen-Cahn equation, namely:
\begin{equation}
    \frac{\partial y}{\partial t} = \nu \frac{\partial^2 y}{\partial x^2} + \mu(y-y^3) + u.
\end{equation}
Here, $\nu>0$ represents again the diffusion coefficient while $\mu$ is the parameter governing the non-linearity. We consider a spatio-temporal domain $\Omega = [a,b]\times[0,T]$ and the following initial and boundary conditions, respectively:
\begin{align}\label{ac}
\begin{aligned}
y(x,0) &= 0.2\sin(\pi x),\qquad && x\in[a,b],\\
y(a,t) &= y(b,t) = 0,\qquad && t\in[0, T].
\end{aligned}
\end{align}
In the cost functional \eqref{cost}, we choose $\bar{y} = y_f = 0$, $\alpha =0.05$, and $\alpha_T=0$, therefore we have 
\begin{equation}
J(y, u) = \frac{1}{2}\int_{0}^{T} \left(  \| y \|_{L^2([a,b])}^2 + 0.05\| u \|^2_{L^2([a,b])} \right) dt \,.
\end{equation}
The adjoint equation is given by
\begin{equation}
    y - \frac{\partial p}{\partial t} = \nu \frac{\partial^2 p}{\partial x^2} +\mu p (1-3y^2),
\end{equation}
equipped with
\begin{align}
\begin{aligned}
p(a,t) &= p(b,t) = 0 , \qquad && t\in[0, T],\\
p(x,T) &= 0,  \qquad && x\in[a,b].
\end{aligned}
\end{align}
Finally, the optimality condition is given by
\begin{equation}\label{eq:opt-ac}
  0.05 u = p.  
\end{equation}
Note that, as in Test 1, it immediately follows that $u$ fulfills the same boundary conditions as $p$. 

As previously, the goal of this numerical test is to solve the control problem and contemporary discover the coefficient $\nu>0$ in \eqref{ac}. We will discuss two scenarios: first, a case in which the parameter to identify is $\nu = 0.1$, for a given $\mu=1$ (case a), and, subsequently, $\nu = 1$, for a given $\mu=11$ (case b). The two test cases reflect two control problems that are deeply different. Specifically, in the second one we enter the region of instability of the solution, which makes the controlled problem considerably more challenging \cite{AGW10}. Furthermore, the spatial domain is given by $a=0$ and $b=1$ in the first case, whereas $a=0$ and $b=2$ in the second case. This choice provides a partially negative initial condition for $y$ in the second case. In both scenarios, the final time is fixed to be $T=0.5$.

This time, we train the OCP-PINN using a scattered dataset of \first{$N_d=12$ (resp. $N_d=14$)} for the uncontrolled state variable \first{(used to solve the inverse problem of finding the unknown parameter)} and \first{$N_d=12$ (resp. $N_d=12$)} training \first{initial condition} points given for the controlled state variable $y$, all selected from the reference solution for Test 2(a) (resp. Test 2(b)). No data is assumed to be available for either the control $u$ or the adjoint $p$, \first{except for the final condition for the latter}.
We evaluate the condition in which we do not have \first{further} information on the initial and terminal conditions of the variables, and we impose boundary conditions at a sparse level, in \first{only $N_b=22$ points for each variable}. Finally, we fix \first{$N_r=520$} residual nodes uniformly distributed in $\Omega$ for evaluating the physical loss terms.
The tolerance for the value of the loss function chosen for the training of this class of test is \first{$\varepsilon_{tol}=10^{-6}$}.

Similarly to Test 1, since the optimality condition \eqref{eq:opt-ac} is given by an algebraic equation, for this test we reduced by one the number of outputs because we imposed the evaluation of $p$ directly as a function of $u$ through \eqref{eq:opt-ac}. This means that, in this setting, the physical residual reads $\mathcal{L}_r = \mathcal{L}_{r,\mathcal{F}} + \mathcal{L}_{r,\mathcal{A}}$, and $\mathcal{O}(p_{NN,r}^n, u_{NN,r}^n, x_r^n, t_r^n;\theta)$, $n=1,\ldots,N_r$, is directly fulfilled. Hence, it is emphasized that, in this test setting, the output variables of the \first{second-in-series PINN} turn out to be only the controlled state and the control. 

\paragraph{Test 2(a): the case $\nu=0.1$, $\mu=1$.}

\begin{figure}[t!]
    \centering    
    \includegraphics[width=0.7\textwidth]{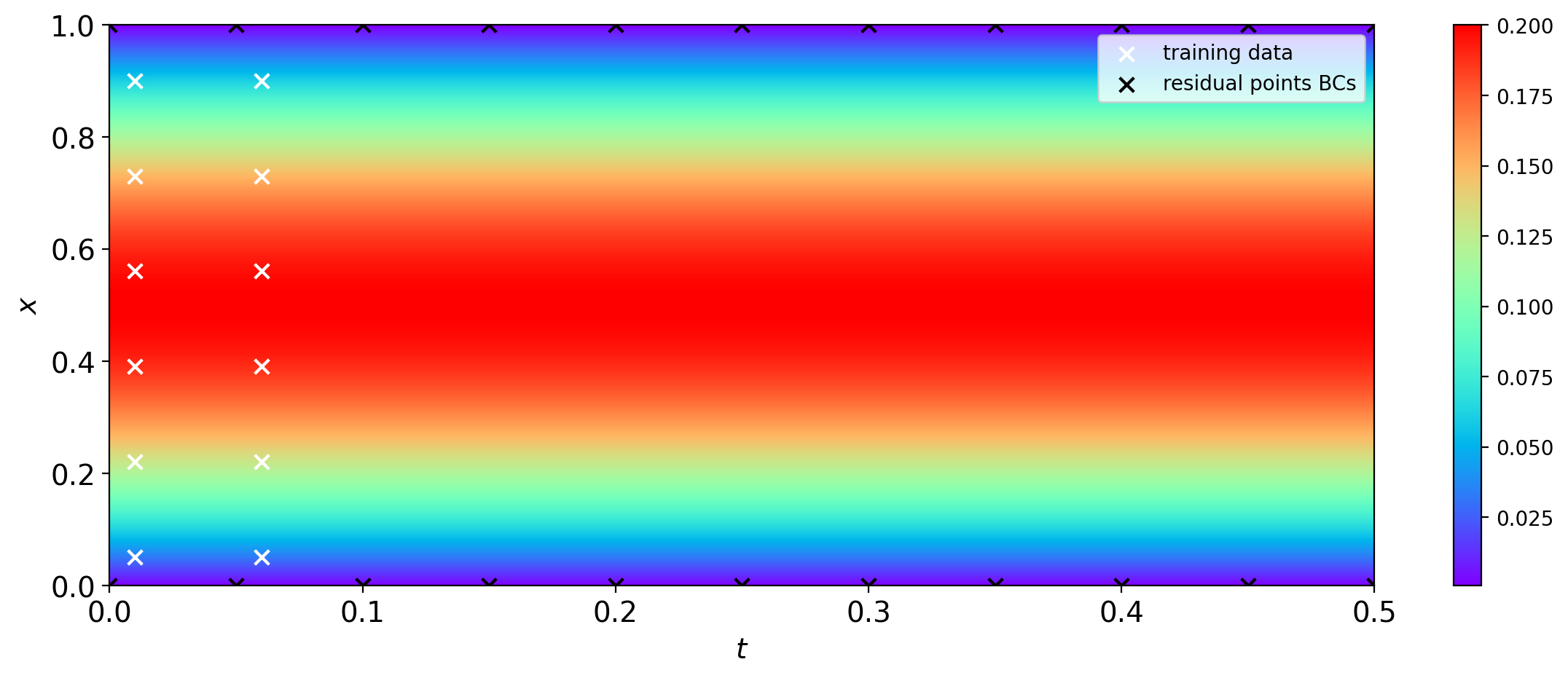}
    \includegraphics[width=0.7\textwidth]{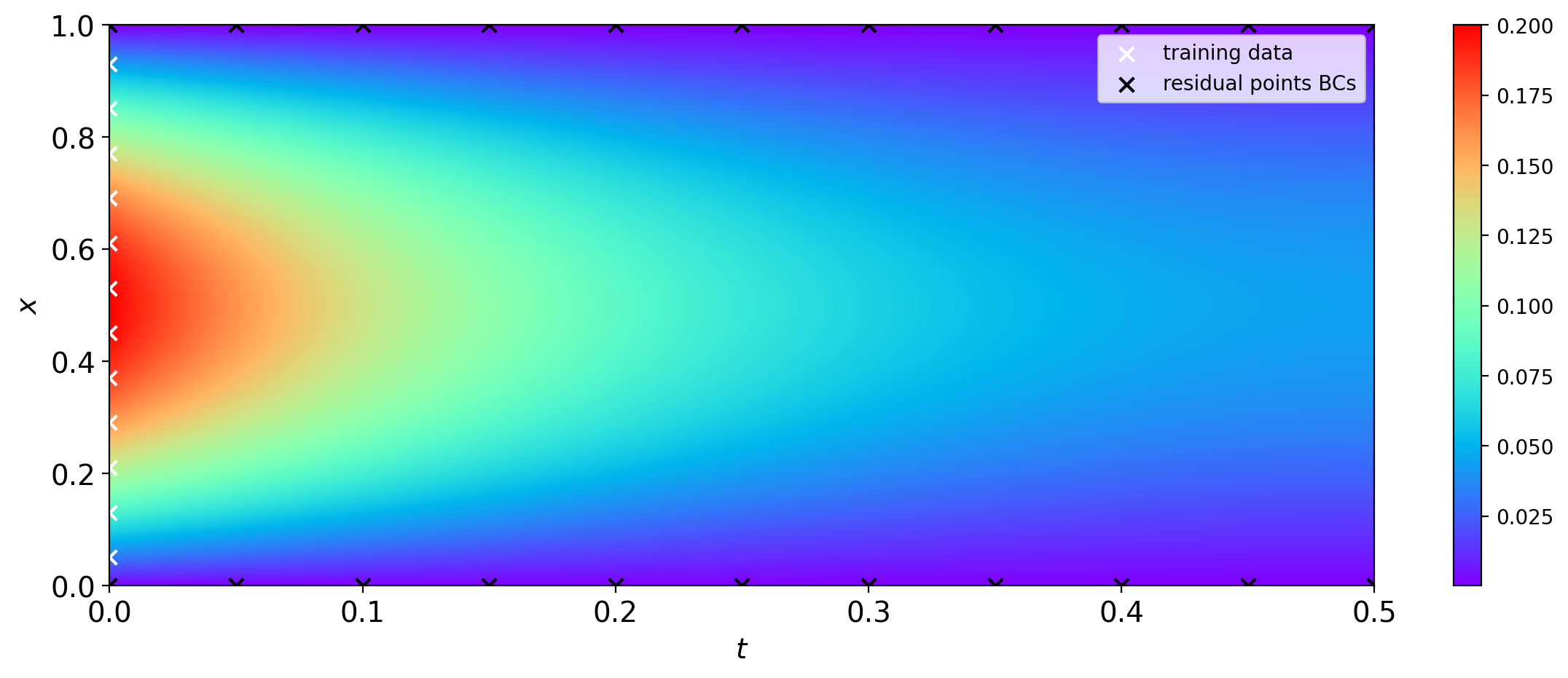}    
    \caption{Test 2(a): Allen-Cahn equation with $\nu=0.1$. Data points of the uncontrolled (top) and controlled (bottom) state variable used for the training of the OCP-PINN.}
    \label{coll-21}
\end{figure}

\begin{figure}[tp!]
    \centering
    \includegraphics[width=0.48\textwidth]{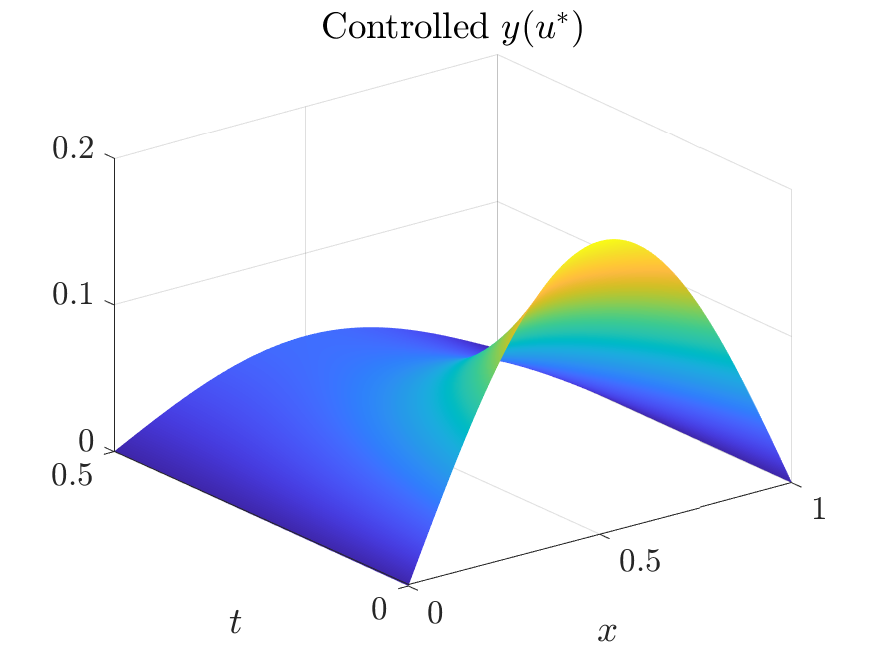}
    \includegraphics[width=0.48\textwidth]{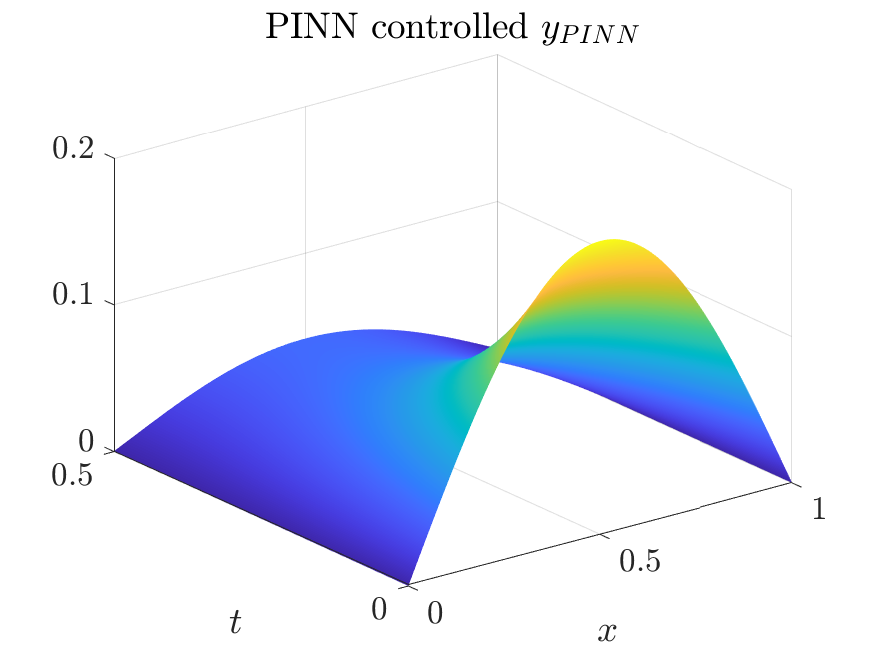}
    \includegraphics[width=0.48\textwidth]{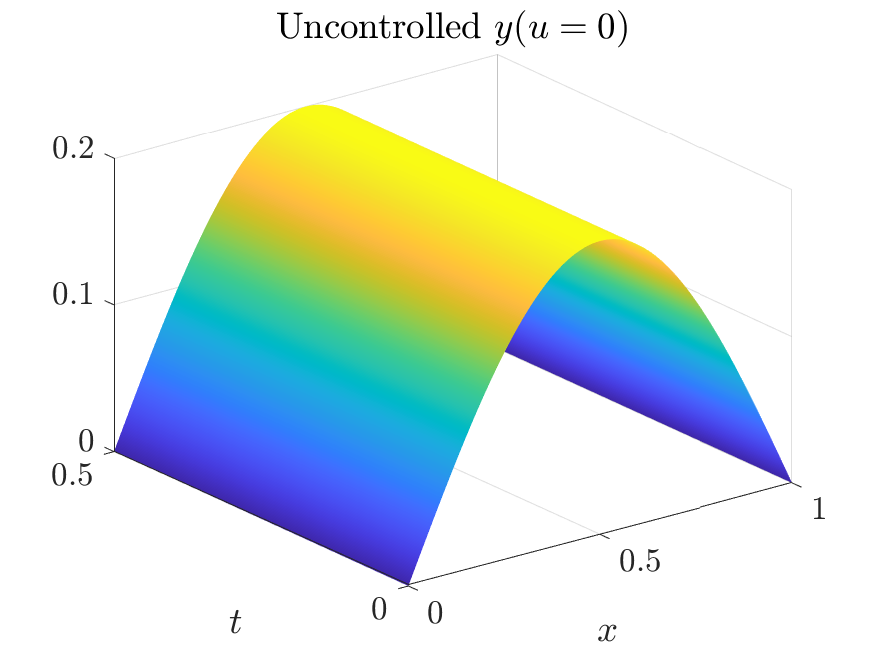}
    \caption{Test 2(a): Allen-Cahn equation with $\nu=0.1$. Reference controlled solution of the state variable (top left), OCP-PINN controlled solution (top right), and reference uncontrolled solution (bottom).}
    \label{fig:ac1}
\end{figure}

\begin{figure}[tp!]
    \centering
    \includegraphics[width=0.48\textwidth]{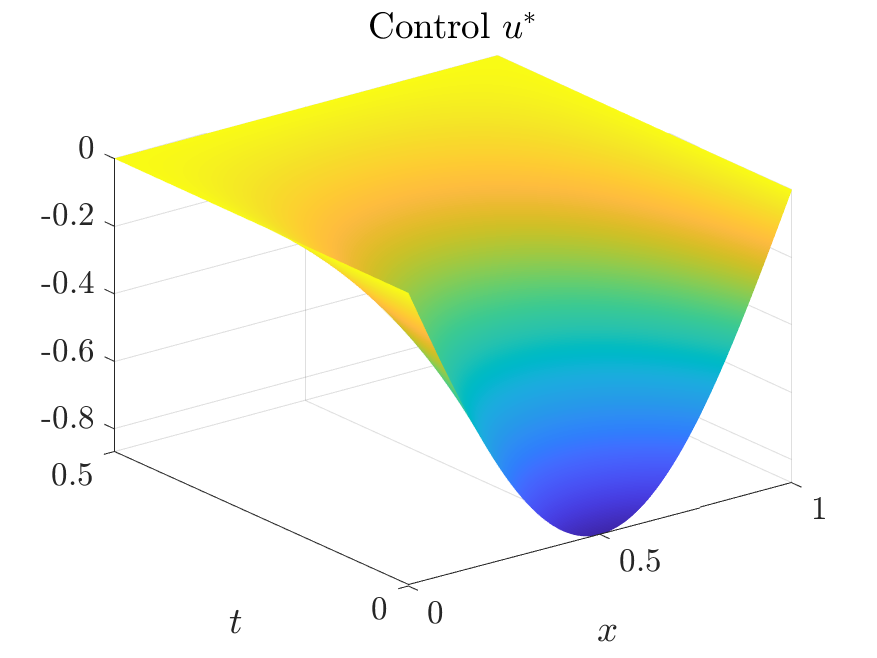}
    \includegraphics[width=0.48\textwidth]{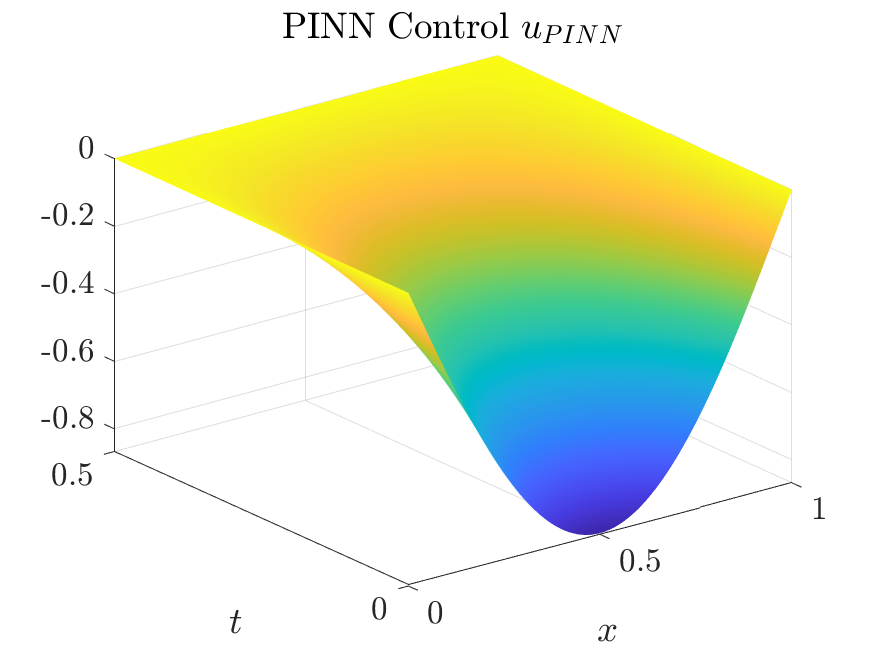}
    \caption{Test 2(a): Allen-Cahn equation with $\nu=0.1$. Reference control variable (left) and control obtained by applying the OCP-PINN (right).}
    \label{ac:contr10}
\end{figure}

\begin{figure}[tp!]
    \centering
    \includegraphics[width=0.48\textwidth]{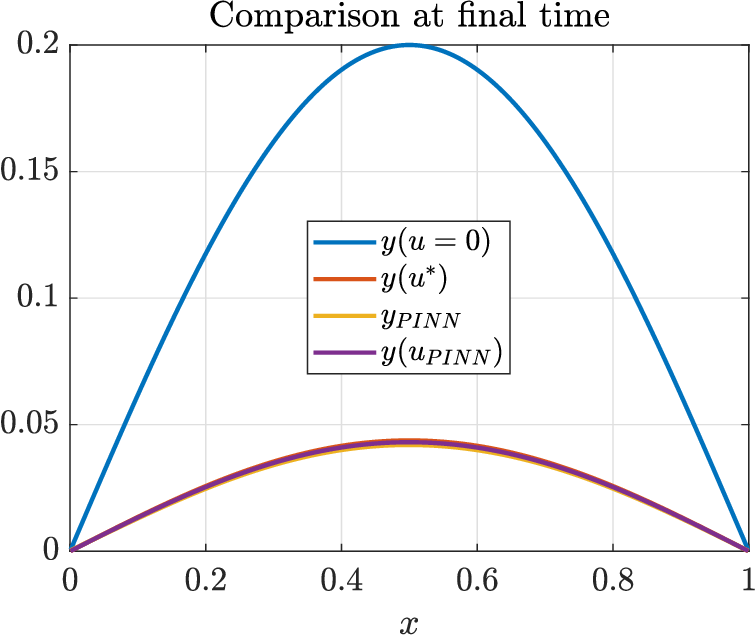}\hspace{0.25cm} 
    \includegraphics[width=0.48\textwidth]{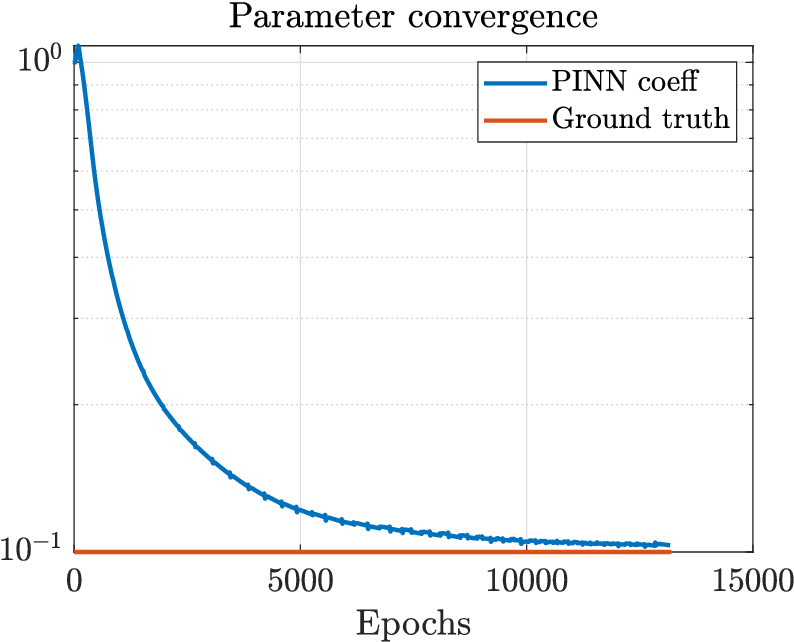}
    \caption{Test 2(a): Allen-Cahn equation with $\nu=0.1$. Comparison of the solutions at final time $T=0.5$ (left) and convergence history for discovering the unknown parameter $\nu$ (right).}
    \label{ac:comp10}
\end{figure}

The dataset used for the training of the \first{two neural networks in series}, for uncontrolled and controlled state variables, is shown in Fig.\ref{coll-21}. 
In Fig. \ref{fig:ac1}, we compare the result obtained with the OCP-PINN for the controlled state with the reference controlled and uncontrolled solutions, while in Fig. \ref{ac:contr10} we present the results in terms of the control variable. We can observe that both the variables $y$ and $u$ are reconstructed very well, apart from a slight reduction in accuracy in identifying the final state of control $u$, which has already been pointed out in the previous test as a small weakness. However, we recall that the control thus returned is identified without having any observed data about it available, only a scattered piece of information at the space boundary.

We compare the solutions of the state variable at the final time $T=0.5$ in the left panel of Fig. \ref{ac:comp10}. Here, it is possible to check the consistent agreement between the reference $y(u^*)$ and the results involving the OCP-PINN method. In the right panel of the same figure, we can appreciate the highly fast convergence of the parameter identification starting with initial guess $\nu_0=1$, which is far from the ground truth. Nonetheless, the discovery of the coefficient performs perfectly, giving \first{$\nu_{PINN} = 0.1032$}. 

The quality of the above-discussed results is confirmed by the relative errors
\first{$\mathcal{E}_1= 0.010697$,  $\mathcal{E}_2=0.003776$, $\mathcal{E}_3= 0.005614$, and $\mathcal{E}_4 = 0.007367$}.
Notice that, once more, all these errors remain of order $O(10^{-2})$. In this simulation, the value of the cost functional computed using the reference solution is $J(y(u^*),u^*)=2.181639$, while for the cost functional obtained with PINN we have $J(y(u_{PINN}),u_{PINN})=2.181596$, with relative error between the two $\mathcal{E}_5 = 1.970995\cdot 10^{-5}$.


\paragraph{Test 2(b): the case $\nu=1$, $\mu=11$.}

\begin{figure}[t!]
    \centering    
    \includegraphics[width=0.8\textwidth]{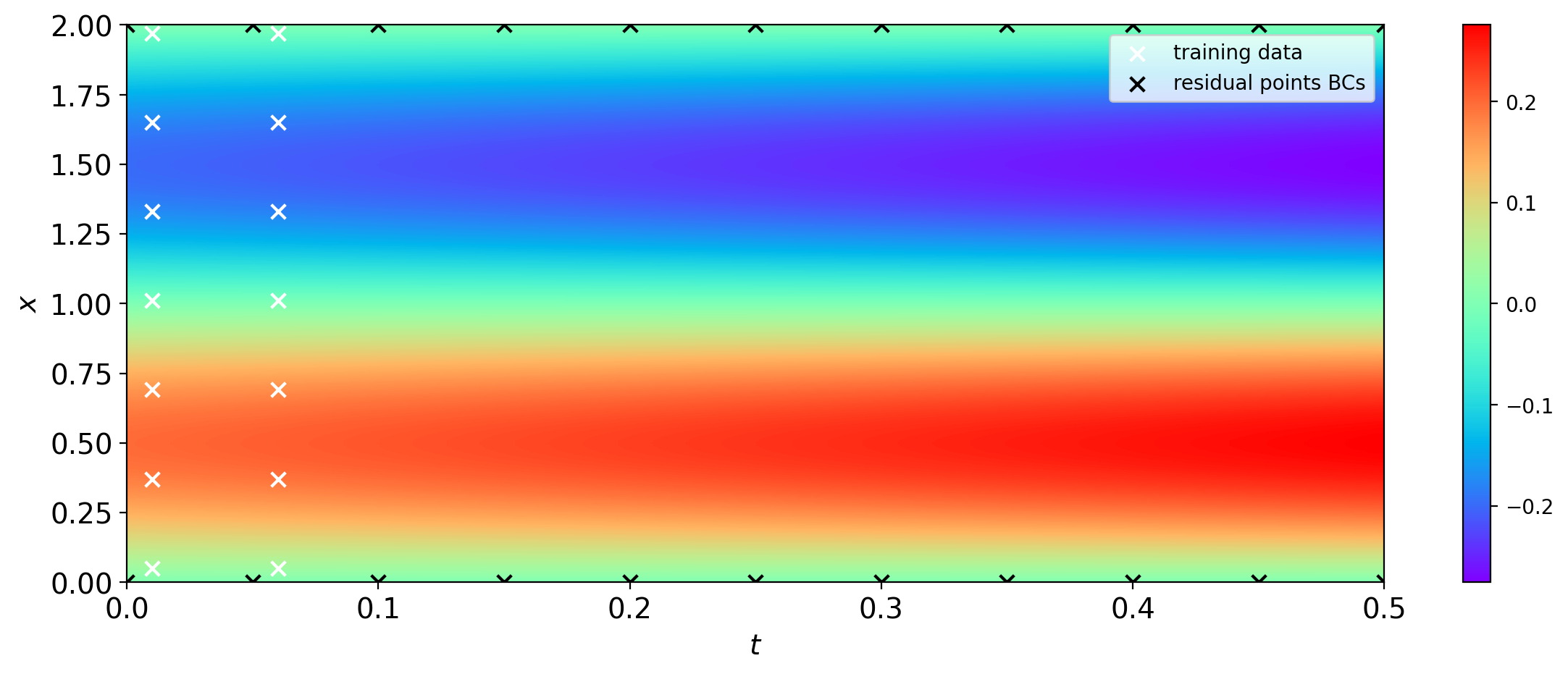}
    \includegraphics[width=0.8\textwidth]{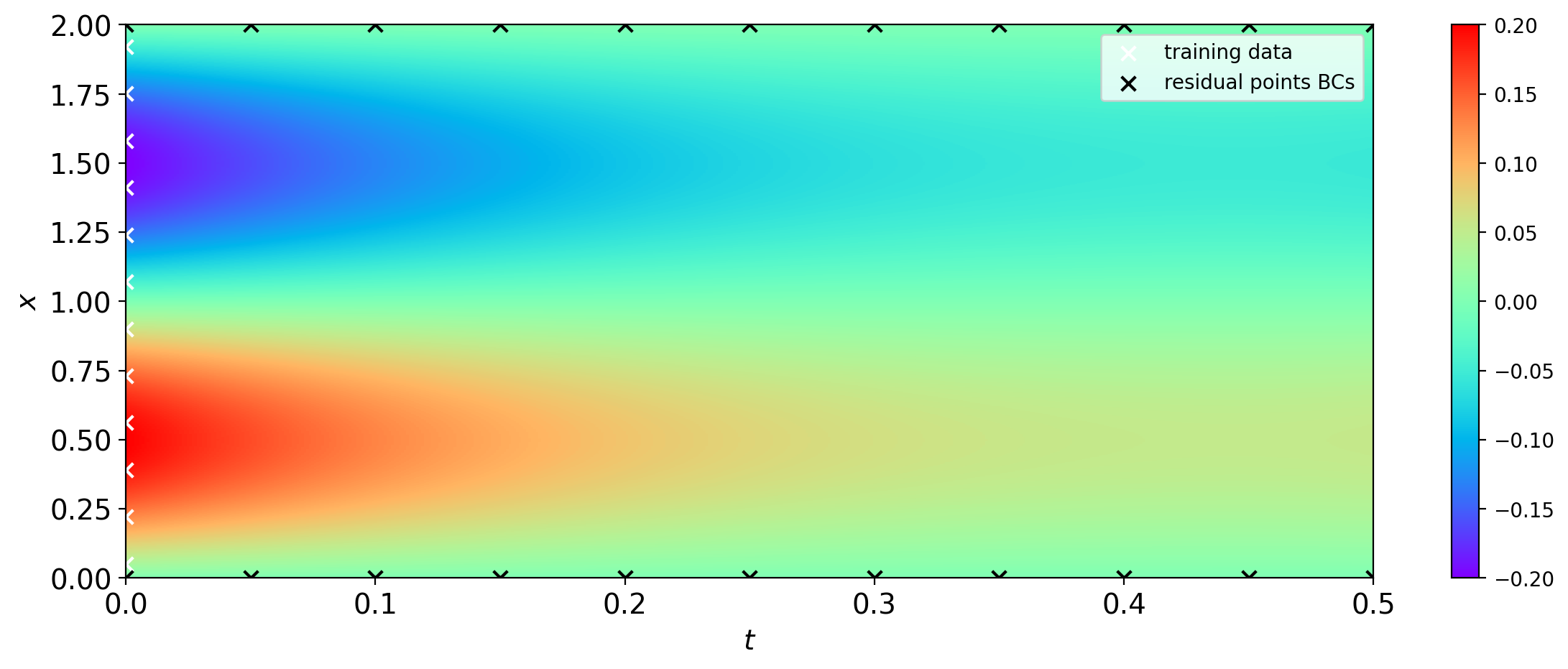}   
    \caption{Test 2(b): Allen-Cahn equation with $\nu=1$. Data points of the uncontrolled (top) and controlled (bottom) state variable used for the training of the OCP-PINN.}
    \label{coll-22}
\end{figure}

\begin{figure}[tp!]
    \centering
    \includegraphics[width=0.48\textwidth]{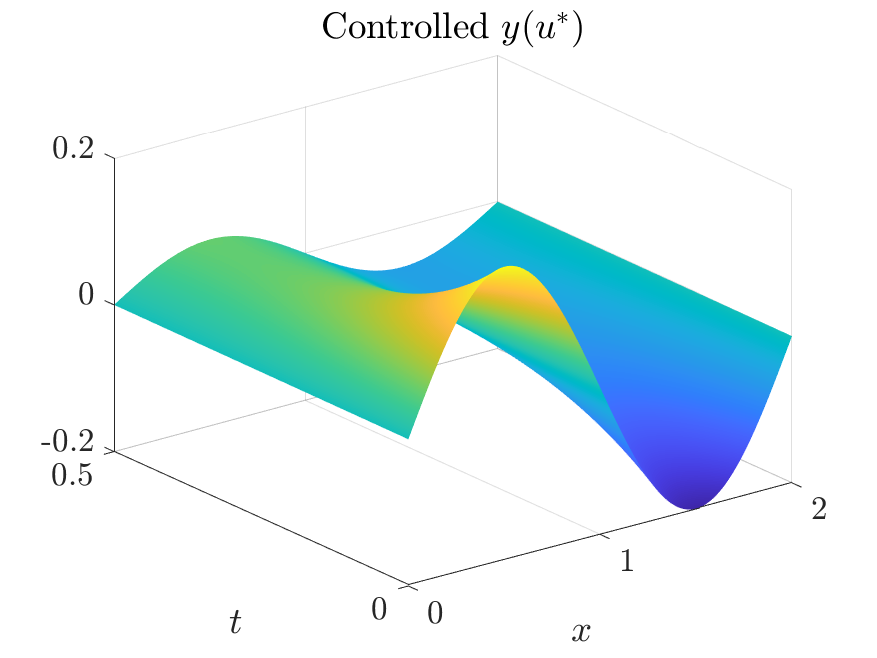}
    \includegraphics[width=0.48\textwidth]{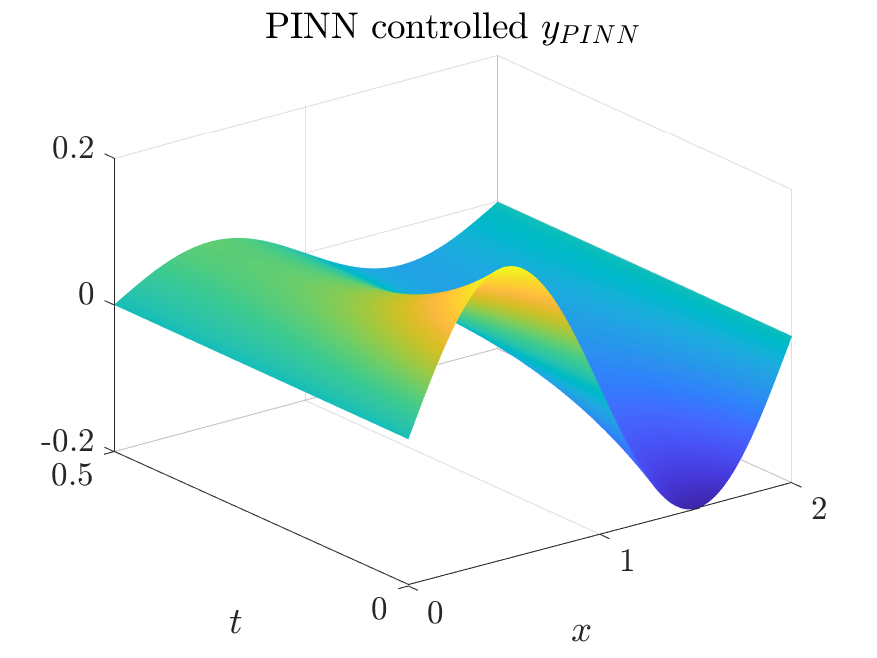}
    \includegraphics[width=0.48\textwidth]{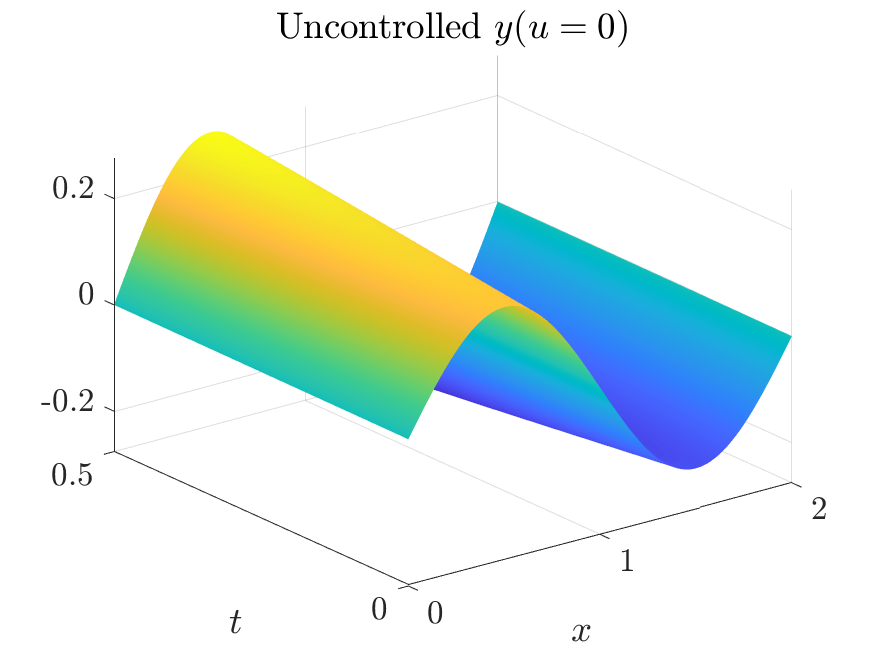}
    \caption{Test 2(b): Allen-Cahn equation with $\nu=1$. Reference controlled solution of the state variable (top left), OCP-PINN controlled solution (top right), and reference uncontrolled solution (bottom).}
    \label{fig:ac2}
\end{figure}

\begin{figure}[tp!]
    \centering
    \includegraphics[width=0.48\textwidth]{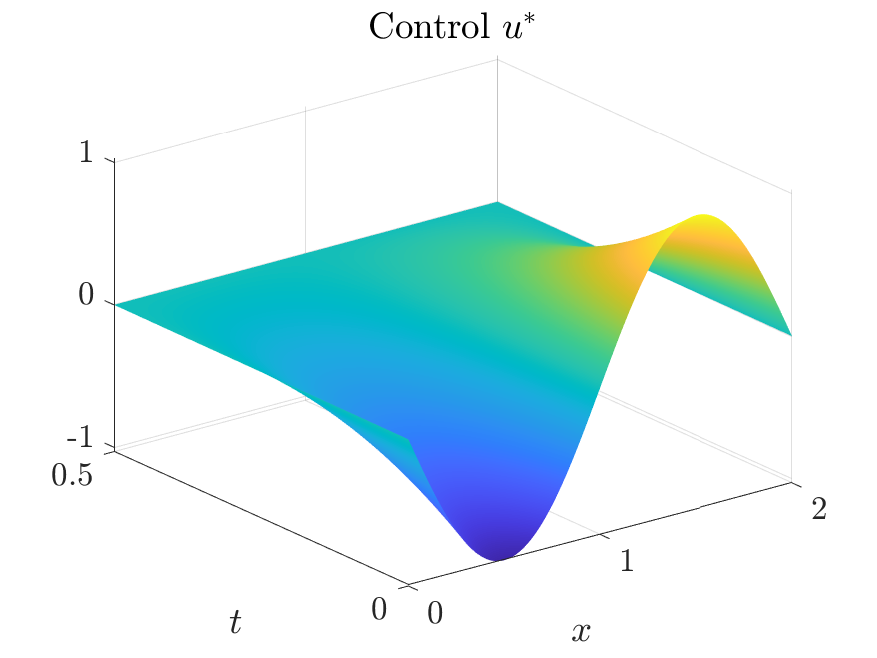}
    \includegraphics[width=0.48\textwidth]{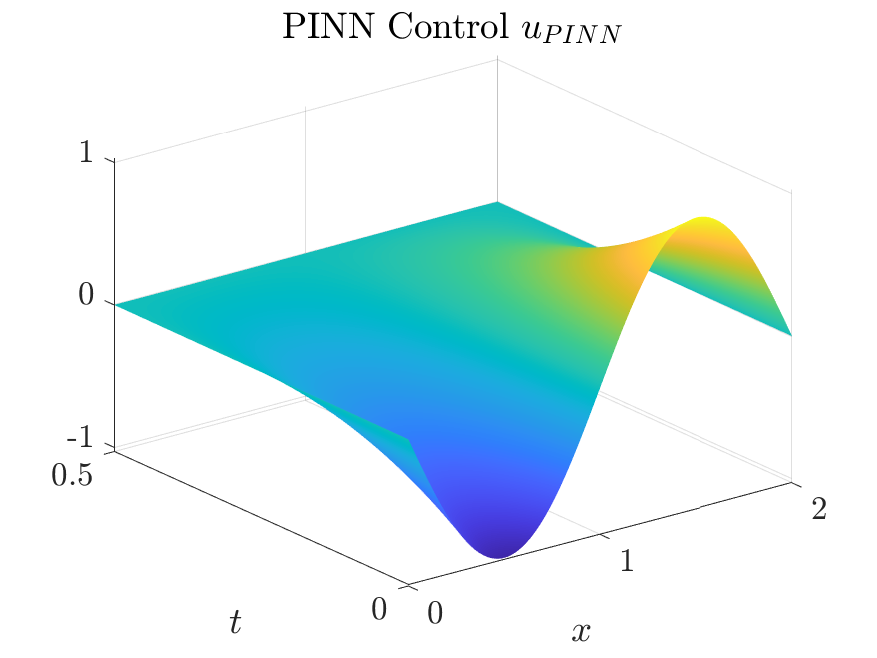}
    \caption{Test 2(b): Allen-Cahn equation with $\nu=1$. Reference control variable (left) and control obtained by applying the OCP-PINN (right).}
    \label{ac:contr1}
\end{figure}

\begin{figure}[tp!]
    \centering
    \includegraphics[width=0.48\textwidth]{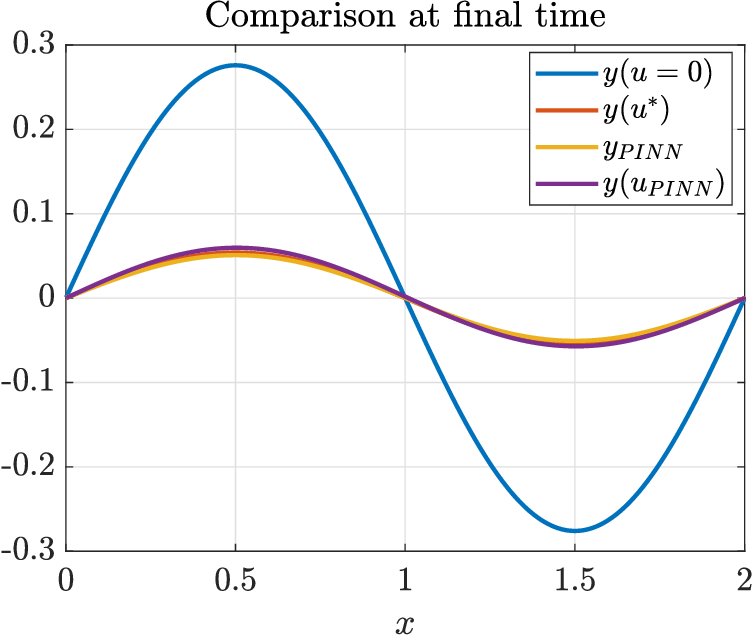}\hspace{0.25cm}
    \includegraphics[width=0.48\textwidth]{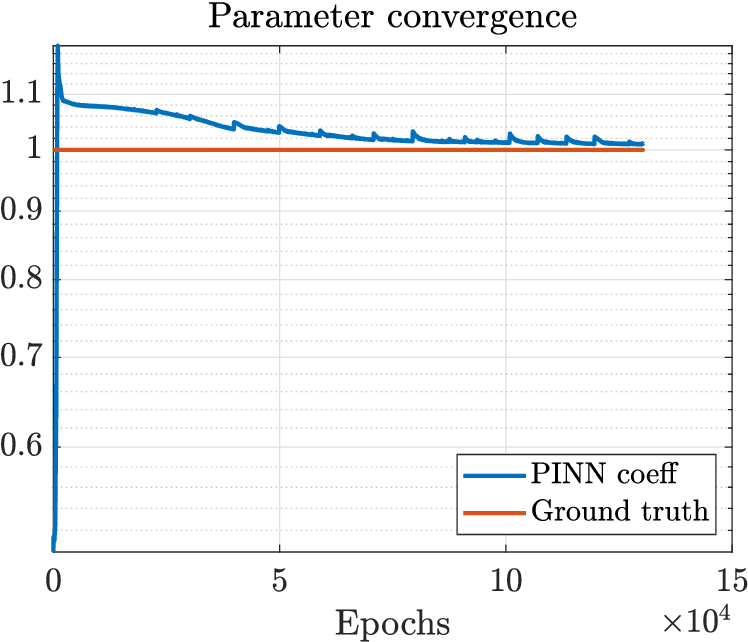}
    \caption{Test 2(b): Allen-Cahn equation with $\nu=1$. Comparison of the solutions at final time $T=0.5$ (left) and convergence history for discovering the unknown parameter $\nu$ (right).}
    \label{ac:comp1}
\end{figure}

We show in Fig. \ref{coll-22} the location of the data points $N_d$ and $N_b$ used to train the OCP-PINN accounting for both uncontrolled and controlled states.
This numerical test concerns a configuration taken from \cite{AGW10} in which the origin is an unstable equilibrium of the problem. Indeed, observing Fig. \ref{fig:ac2} it can be noticed that the amplitude of the uncontrolled solution is increasing over time and the control acts to reverse this trend by bringing the variable to decrease, which makes the control problem extremely complicated. Nevertheless, from Figs. \ref{fig:ac2} and \ref{ac:contr1} it can be appreciated the remarkable agreement the OCP-PINN solution presents with the reference one.


In Fig. \ref{ac:comp1}, it can be seen that the final profile of the controlled solution matches, with a very small error, the reference curve, and the rapid convergence of the identification of the unknown parameter. Indeed, in this test, the parameter discovery turns out to be \first{ as fast as in} the previous case because the neural network identifies the parameter \first{$\nu_{PINN}=1.0100$} with great accuracy already after \first{$1.3\cdot 10^4$} epochs (initial guess $\nu_0=0.5$). 

Finally, we report the following relative norms:
\first{$\mathcal{E}_1= 0.016978 $, $\mathcal{E}_2= 0.031128 $, $\mathcal{E}_3 =0.017843$, $\mathcal{E}_4= 0.046732$}, which remain of order $O(10^{-2})$. Concerning the cost functional, its value computed using the reference solution is $J(y(u^*),u^*)=5.425312$, while for the one obtained with PINN we have $J(y(u_{PINN}),u_{PINN})=5.429540$. The relative error is $\mathcal{E}_5 = 7.793100\cdot 10^{-4}$.



\subsection{Test 3: Korteweg--de Vries equation}\label{sec:test3}
The last numerical simulation concerns the control of the Korteweg--de Vries (KdV) equation  with an additional second-order derivative term that acts as a supplementary diffusion \cite{Kato1979,WAZWAZ2008485}:
\begin{equation}\label{eq:kdv}
\frac{\partial y}{\partial t} + \mu y\frac{\partial y}{\partial x} = \nu \frac{\partial^2 y}{\partial x^2} - \lambda \frac{\partial^3 y}{\partial x^3} + u.
\end{equation}
This equation models the behavior of small dispersive waves: $y$ identifies the wave profile, the drift term accounts for the non-linearity in the wave propagation (indicating that the speed of the wave depends on its amplitude), the second-order partial derivative with respect to $x$ describes a diffusion effect, while the third-order partial derivative represents the dispersion. The balance between the nonlinear drift term and the diffusive and dispersive ones is what allows the formation of stable solitary waves, or \textit{solitons}, which is a key feature of the KdV equation.  We focus on the spatio-temporal domain $\Omega = [a,b] \times [0,T]$, with $a=-4$, $b=4$, $T=0.5$, and set the following initial and boundary conditions for the state variable:
\begin{align}\label{kdvbc}
\begin{aligned}
y(x,0) &= \sqrt{ \frac{2}{\pi}}e^{-2x^2},\qquad && x\in[a,b],\\
y(a,t) &= y(b,t) = 0,\qquad && t\in[0, T].
\end{aligned}
\end{align}
In the cost functional \eqref{cost}, we set $\bar{y} = y_f = 0$, $ \alpha=0.1$, and $\alpha_T=0$, so that we have 
\begin{equation}
J(y, u) = \frac{1}{2}\int_{0}^{T} \left(  \| y \|_{L^2([a,b])}^2 + 0.1\| u \|^2_{L^2([a,b])} \right) dt\,.
\end{equation}
The adjoint equation results
\begin{equation}
    y - \frac{\partial p}{\partial t} - \mu y \frac{\partial p}{\partial x} = \nu \frac{\partial^2 p}{\partial x^2} + \lambda \frac{\partial^3 p}{\partial x^3},
\end{equation}
coupled with
\begin{align}\label{eq:kdvadj}
\begin{aligned}
p(a,t) &= p(b,t) = 0 , \qquad && t\in[0, T],\\
p(x,T) &= 0,  \qquad && x\in[a,b].
\end{aligned}
\end{align}
Finally, we have the optimality condition
\begin{equation}\label{eq:opt-kdv}
  0.1u = p.  
\end{equation}
The chosen values of the parameters appearing in \eqref{eq:kdv} and \eqref{eq:kdvadj} are $\mu = -5$,  $\nu = 1$, and $\lambda = 1$.

\begin{figure}[tp!]
    \centering    
    \includegraphics[width=0.7\textwidth]{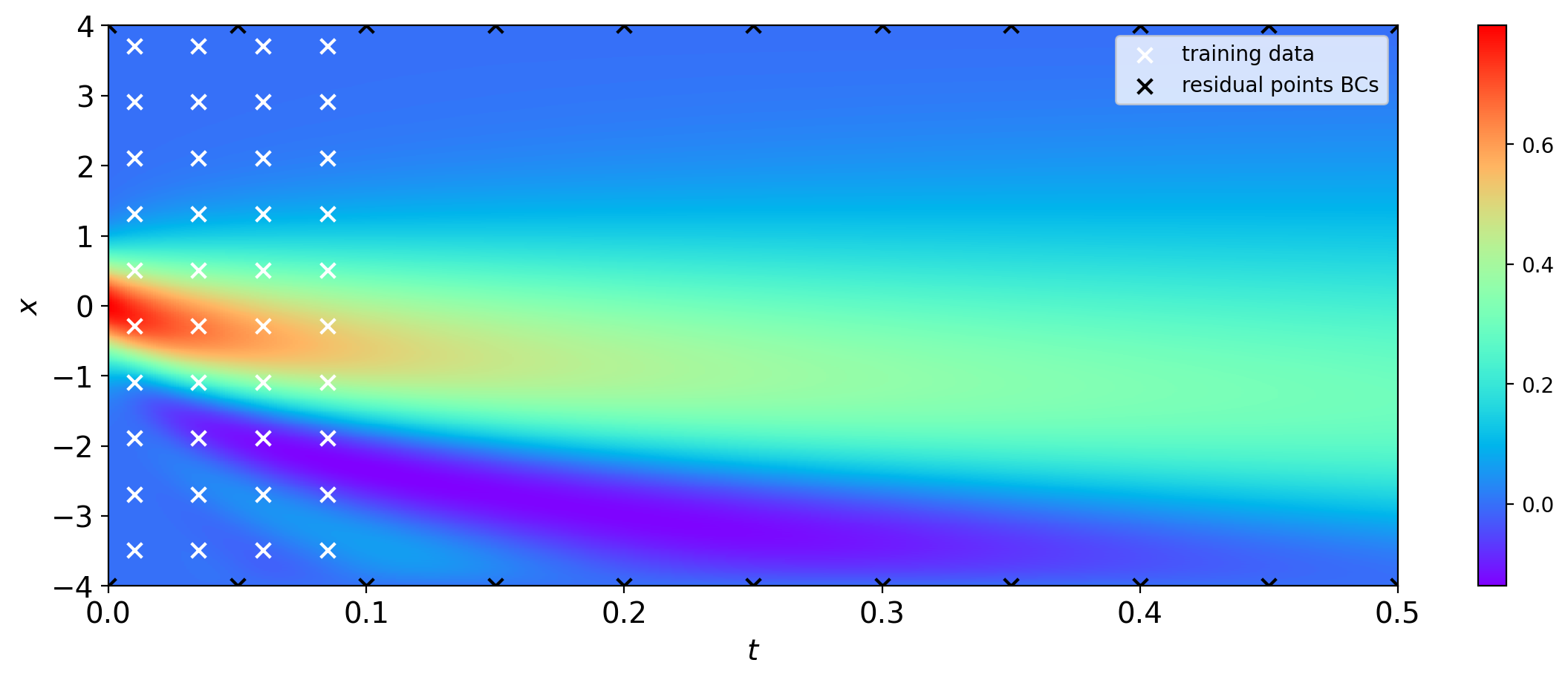}
    \includegraphics[width=0.7\textwidth]{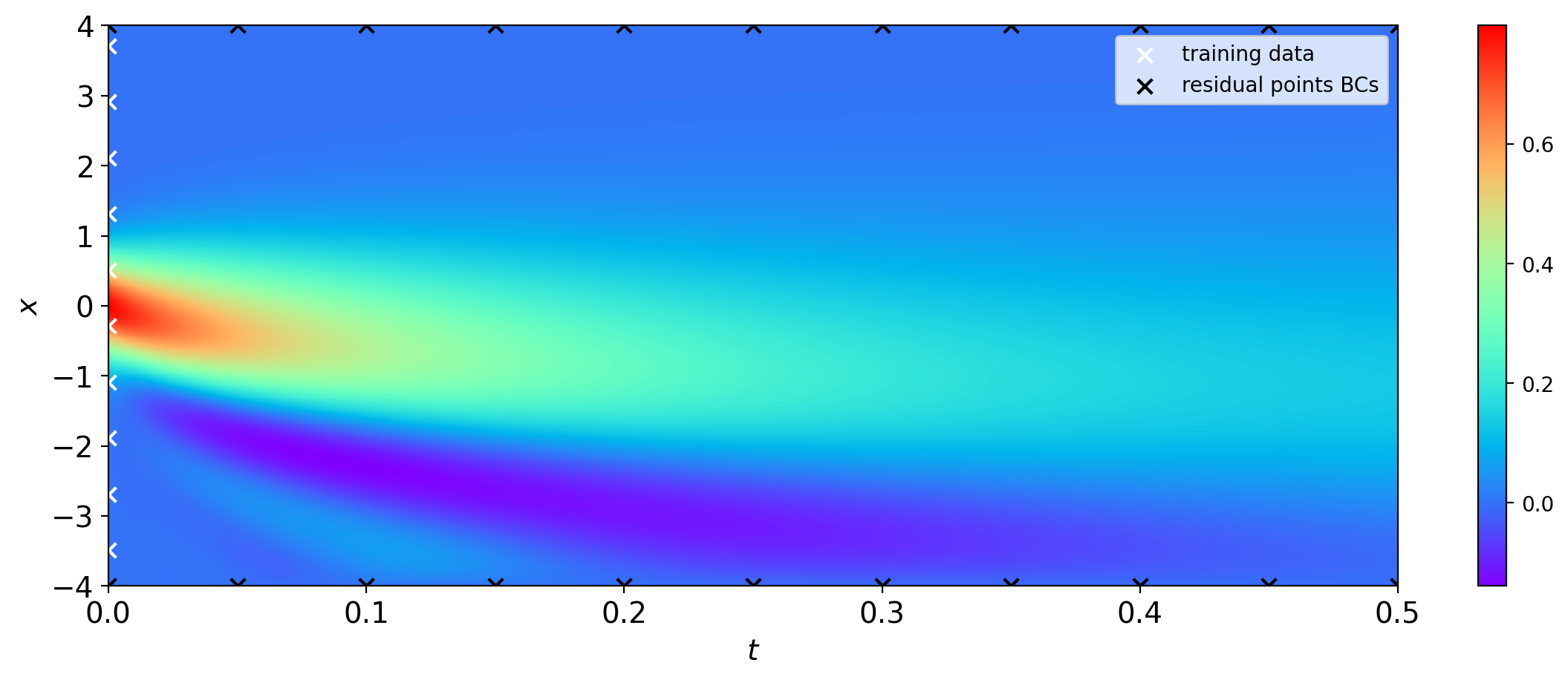}   
    \caption{Test 3: Korteweg--de Vries equation. Data points of the \first{uncontrolled (top) and controlled (bottom)} state variable used for the training of the OCP-PINN.}
    \label{fig:kdv_data}
\end{figure}

\begin{figure}[tp!]
    \centering
    \includegraphics[width=0.48\textwidth]{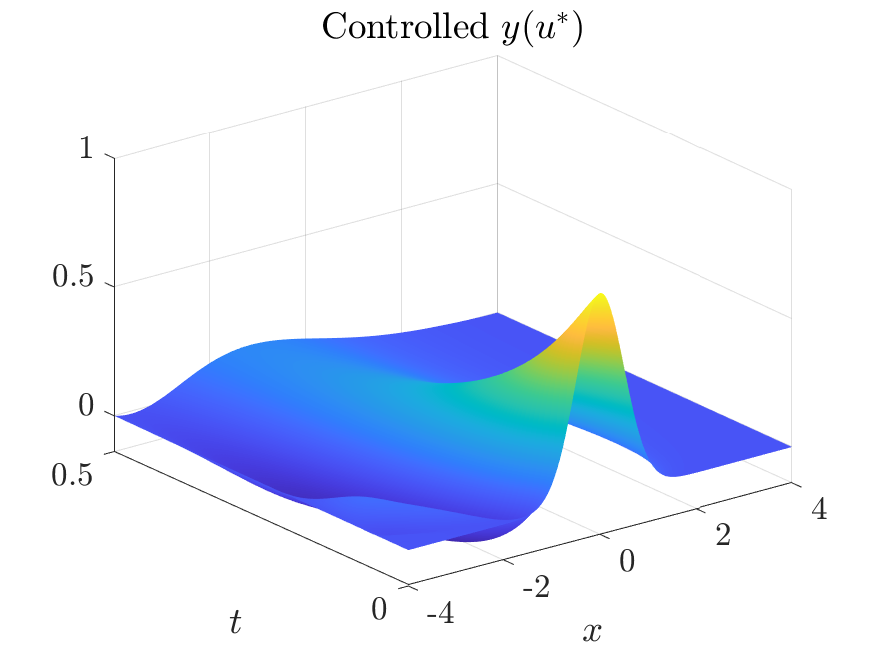}
    \includegraphics[width=0.48\textwidth]{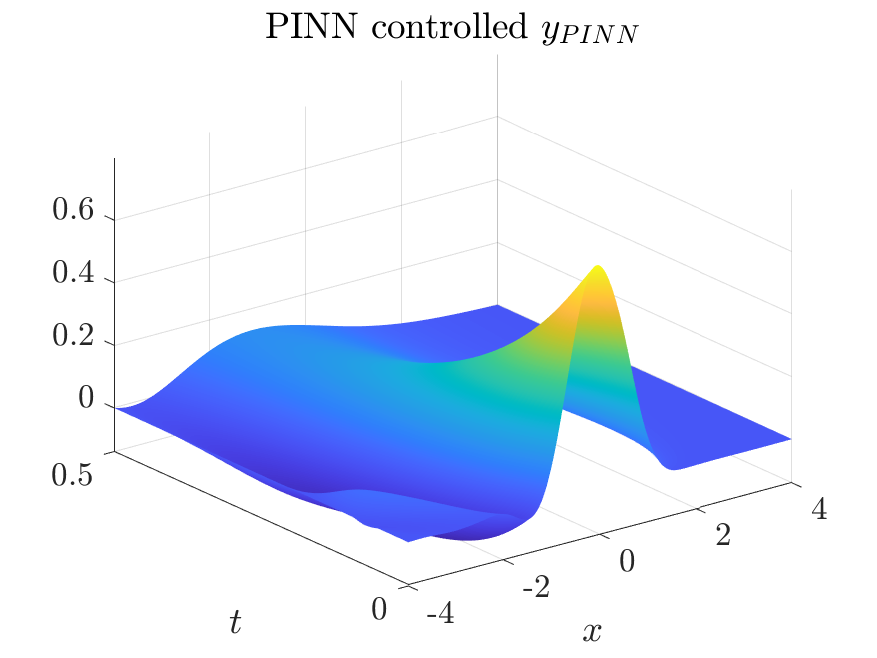}
    \includegraphics[width=0.48\textwidth]{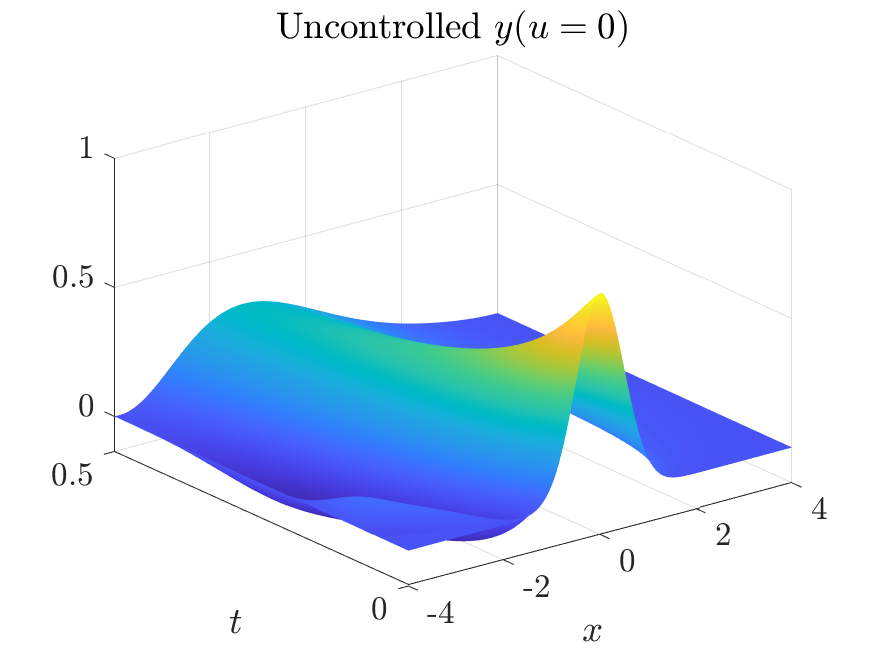}
    \caption{Test 3: Korteweg--de Vries equation. Reference controlled solution of the state variable (top left), OCP-PINN controlled solution (top right), and reference uncontrolled solution (bottom).}
    \label{fig:kdv1}
\end{figure}

\begin{figure}[tp!]
    \centering
    \includegraphics[width=0.48\textwidth]{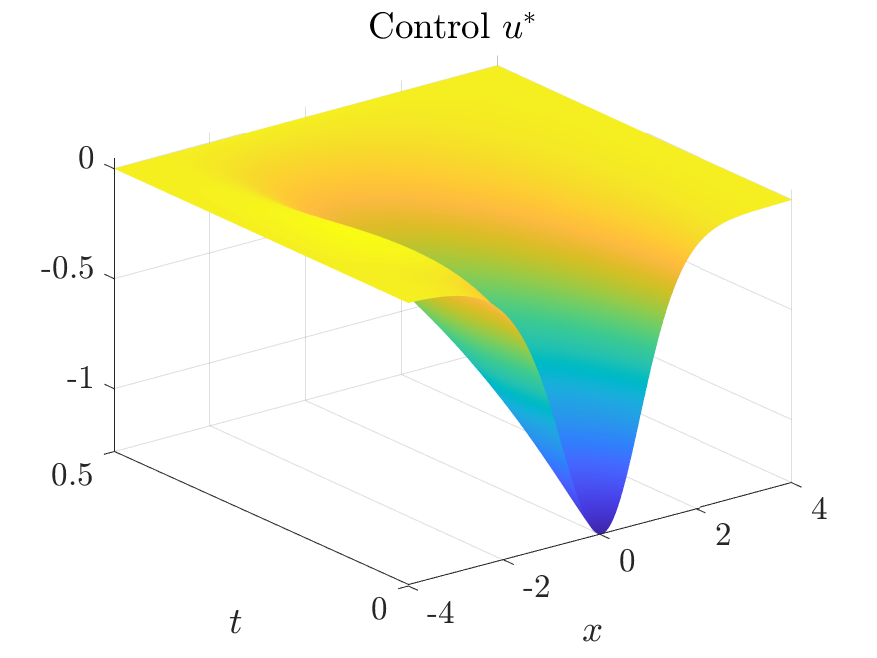}
    \includegraphics[width=0.48\textwidth]{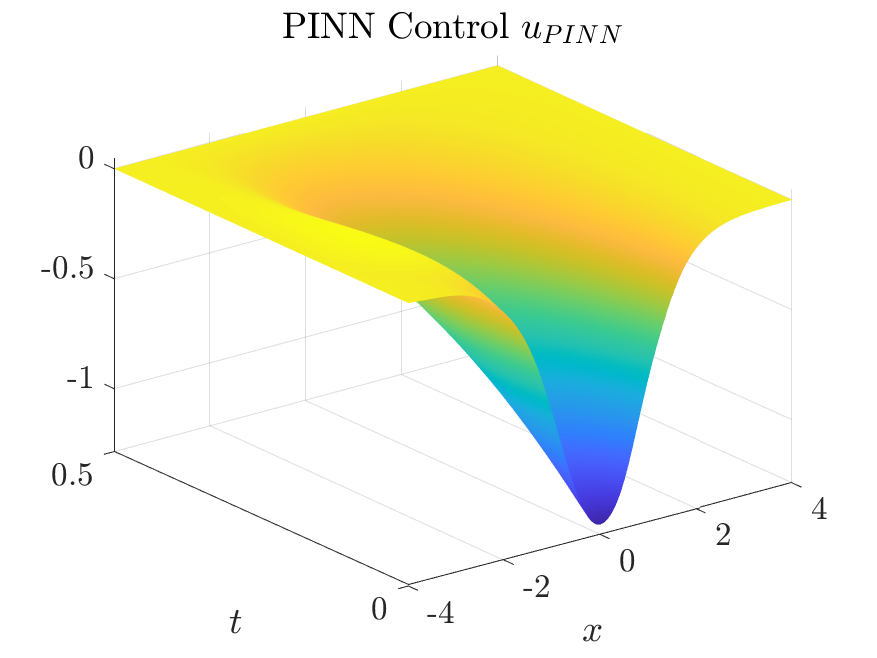}
    \caption{Test 3: Korteweg--de Vries equation. Reference control variable (left) and control obtained by applying the OCP-PINN (right).}
    \label{fig:kdvcontr}
\end{figure}

\begin{figure}[tp!]
    \centering
    \includegraphics[width=0.48\textwidth]{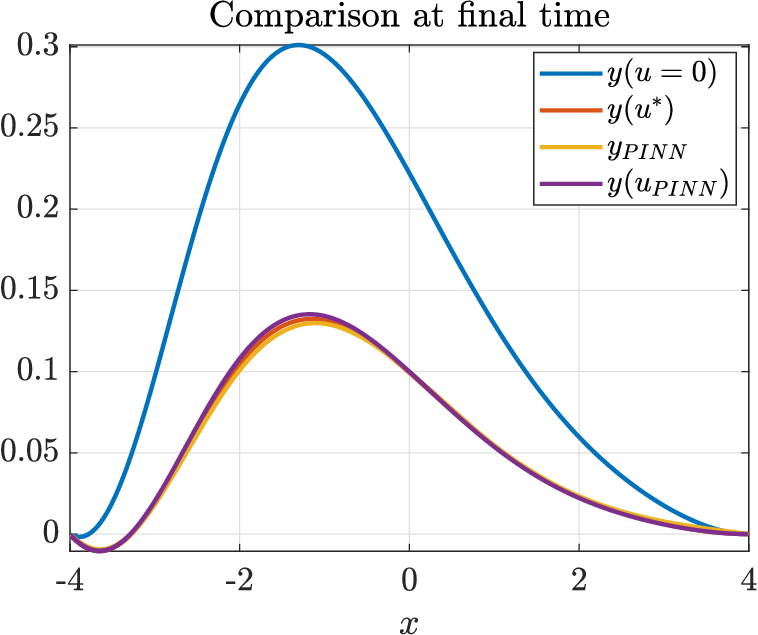}\hspace{0.25cm}
    \includegraphics[width=0.48\textwidth]{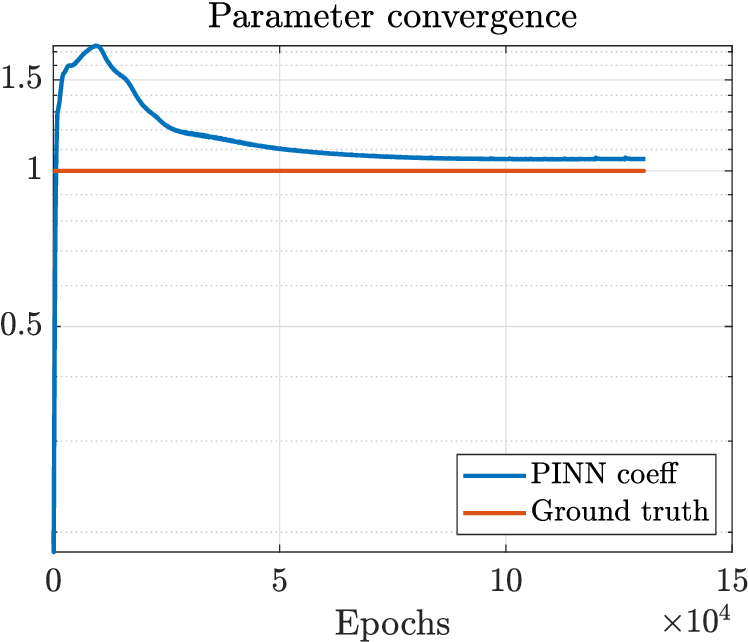}
    \caption{Test 3: Korteweg--de Vries equation. Comparison of the solutions at final time $T=0.5$ (left) and convergence history for discovering the unknown parameter $\nu$ (right).}
    \label{kdv:comp}
\end{figure}

We want again to solve the control problem and, at the same time, learn the value of the diffusion coefficient $\nu > 0$, which is considered unknown. 
As in the previous cases, the optimality condition \eqref{eq:opt-kdv} is given by an algebraic equation. Therefore, for this test, we reduce the number of outputs of the \first{second-in-series PINN} to 2 and directly evaluated $u$ as a function of $p$ through Eq. \eqref{eq:opt-kdv}, as in Test 2. 

We train the OCP-PINN using a scattered dataset of \first{$N_d=40$ training points for the uncontrolled state $y_{unc}$ and only $N_d=10$} training points given for the \first{initial state of the} controlled variable $y$, selected from the reference solutions, while no data is considered available for either the control $u$ or the adjoint $p$. We impose boundary conditions at a sparse level, in \first{$N_b=22$} points, and for the adjoint variable we consider available terminal conditions. The location of training points is shown in Fig. \ref{fig:kdv_data}. Finally, we fix \first{$N_r=520$} residual nodes uniformly distributed in $\Omega$ for evaluating the physical loss terms.
Training is performed for a number of iterations (epochs) such that an error tolerance of \first{$\varepsilon_{tol}=10^{-6}$} is satisfied by the loss function.

Fig. \ref{fig:kdv1} shows a comparison between the reference controlled solution, the OCP-PINN solution, and the reference uncontrolled state variable, while Fig. \ref{fig:kdvcontr} presents a comparison of the results in terms of the control variable. We can observe, also in this test, a good agreement between the reference solutions and the variables reconstructed with the proposed numerical technique.

We compare the solutions of the state variable $y$ at the final time $T=0.5$ in the left panel of Fig. \ref{kdv:comp}, again showing a remarkable agreement among the controlled states computed with the different approaches. On the right panel of the same Figure, we show the convergence of the discovered diffusion coefficient. The final predicted value is \first{$\nu_{PINN} = 1.0547$}, which is satisfactorily close to the ground truth $\nu = 1$ (initial guess $\nu_0 = 0.2$). We can also notice that, once more, the convergence of the parameter does not present evident spurious oscillations.

The relative errors remain of order $O(10^{-2})$ and are equal to \first{$\mathcal{E}_1 =  0.038560$,  $\mathcal{E}_2=0.010061$, $\mathcal{E}_3= 0.020784$, and $\mathcal{E}_4 = 0.038614$}. The value of the cost functional computed with the reference solution is $J(y(u^*),u^*)=159.488167$, while for the cost functional obtained using PINN we have $J(y(u_{PINN}),u_{PINN})=154.386241$, with relative error $\mathcal{E}_5 = 3.198937\cdot 10^{-2}$.


\section{Conclusions}
\label{sec:final}
Physics-Informed Neural Networks have emerged in recent years as a powerful tool for solving differential equations by integrating data and physical laws into the neural network training process.
In this paper, we propose an extension of the standard PINN framework, exploiting the Lagrange multipliers approach, to allows its use for solving optimal control problems under conditions of partial physical information about the investigated dynamics. Specifically, the OCP-PINNs proposed here allow the identification of the control variable and any unknown model parameters with an online process.

The validity of the proposed methodology has been confirmed by several numerical tests.
However, we would like to take advantage of this concluding section to report on various aspects and considerations that emerged during this study regarding the application of PINNs to control problems, discussing their advantages, challenges, and potential future applications, hoping that this will better direct any readers interested in undertaking the exploration of these innovative techniques. 
We believe that many of the issues reported below represent shared doubts in the numerical community and certainly deserve special attention. In this context, we limit ourselves to raising a few brief comments on them, knowing that the discussion here presented cannot be considered exhaustive.

\begin{itemize}
    \item \textit{When is it better to use PINNs instead of traditional numerical methods?} 
    This is certainly a difficult question to answer given the complexity of the possible applications involved in using these techniques. However, we would like to say, without claiming in any way to have in our hands the definitive answer to this question, that PINNs prove to be a useful tool when one wants to solve high-dimensional problems, complex geometries, or scenarios where some information about the problem, such as initial and boundary conditions, is missing and it is therefore impracticable the application of traditional numerical methods (such as Finite Differences, Finite Volume and Finite Elements).\\

    \item \textit{To train a PINN, do you need data, and if so, how much and what kind?} 
    While PINNs can work with limited data, the quality and quantity of data can significantly impact the accuracy of the results. Typically, PINNs require scattered initial or boundary conditions and some additional data points to guide the training. The amount needed depends on the problem's complexity and the quality of the prior physical information encoded in the model. Not surprisingly, it has been observed that the more data provided, the faster the convergence to the desired result.\\

    \item \textit{Well-posedness of the inverse problem: How can we ensure that the inverse problem has a unique and stable solution?}
    Ensuring well-posedness involves verifying that the inverse problem satisfies existence, uniqueness, and stability criteria. Regularization techniques and careful formulation of the loss function can help stabilize the solution. Additionally, leveraging prior knowledge and constraints can improve the well-posedness.\\

    \item \textit{How much problem-dependent is the choice and design of the neural network?}
    We regret to say that the design of the neural network, including its architecture, activation functions, and training procedures, appears to depend heavily on the specific problem. Factors such as the type of differential equations (e.g., elliptic, parabolic), the geometry of the domain, and the available data influence the design of the network, and customization of elements such as size and learning rate for the problem at hand is critical to achieving good performance. This makes each problem different and requires, frequently, a trial-and-error phase of configuration that we understand is certainly not desirable.  It is worth reporting here that, at the early stage of this work, we started to study the control of the linear heat equation by OCP-PINNs. During this preliminary analysis, we could not clearly identify any kind of monotone convergence with respect to the number of nodes or hidden layers used in the network setting.\\

    \item \textit{The importance of mathematical foundations and physical principles in the design and training of neural networks.}
    The incorporation of mathematical foundations and physical principles ensures that PINNs respect the underlying laws governing the system, allowing for improved interpretability, accuracy, and generalization of the model. However, it is important to keep in mind that the formulation of the residual terms in the loss function must be done accurately to ensure that they can guide the neural network toward physically plausible solutions consistent with known behaviors (see, for example, \cite{Bertaglia2022b,Liu2024,Jin2024} for the formulation of asymptotic-preserving neural networks for multiscale problems). \\

    \item \textit{What is the acceptable order of error for solutions obtained using PINNs compared to traditional methods?}
    The error given by PINNs, relative to an exact reference solution, typically ranges from \(O(10^{-2})\) to \(O(10^{-3})\), depending on the requirements of the problem. Therefore, PINNs turn out to yield errors that are not exactly low. While traditional numerical methods can achieve higher accuracy for well-posed and well-understood problems, it should still be kept in mind that PINNs offer flexibility and robustness in handling complex scenarios where such accuracy might be difficult to achieve anyway.\\

    \item \textit{Challenges in accurately capturing boundary conditions and physical constraints.}
    Accurately capturing boundary conditions and physical constraints can be challenging because of the complexity of ensuring strict adherence to the physics of the neural network especially with very little information available. In particular, we found considerable difficulty for neural networks in capturing solutions very close to zero, a difficulty that is likely explained by the orders of error achievable by PINNs discussed in the previous section. Another well-known weakness of neural networks is that they tend to return softer solutions than expected, which makes them not very attractive for solving hyperbolic problems with shock wave formation, for example.\\
    
    \item \textit{What is the expected computational cost of using neural networks?}
    It cannot be hidden that one of the main disadvantages of neural networks, compared to conventional numerical methods, is the computational cost of the training phase, which is generally very high. The advantage of using PINNs lies in the possibility of re-using the trained models after carrying out a deep and extensive training process, which results in a network robust enough to be re-applied to different studies. That said, in case histories where solving the system of equations is already very expensive with classical methods (think for example of the case of high-dimensional kinetic equations \cite{Jin2021}), we can say that neural networks can still be competitive in terms of computational efficiency. Not to mention cases where the application of traditional numerical methods is not feasible due to a lack of physical knowledge about the dynamics to be studied, as previously discussed. 
\end{itemize}

\section*{Acknowledgements}
The authors have developed this work within the activities of the project ``Data-driven discovery and control of multi-scale interacting artificial agent systems” (code P2022JC95T), funded by the European Union -- NextGenerationEU, National Recovery and Resilience Plan (PNRR) – Mission 4 component 2, investment 1.1 ``Fondo per il Programma Nazionale di Ricerca e Progetti di Rilevante Interesse Nazionale" (PRIN).
The support of the INdAM--GNCS group is also acknowledged. A.A. is part of the GNCS Project ``Metodi numerici innovativi per equazioni di Hamilton-Jacobi” (code E53C23001670001); G.B. and E.C. participate in the GNCS Project ``Metodi numerici per le dinamiche incerte" (code E53C23001670001).
G.B. has also been funded by the European Union -- NextGenerationEU under the program ``Future Artificial Intelligence -- FAIR" (code PE0000013), MUR PNRR, Project ``Advanced MATHematical methods for Artificial Intelligence -- MATH4AI".
This manuscript reflects only the authors’ views and opinions, neither the European Union nor the European Commission can be considered responsible for them.
A.A. is also supported by MIUR with PRIN
project 2022 funds (P2022238YY5, entitled ‘‘Optimal control problems: analysis, approximation’’)

\bibliographystyle{abbrv}
\bibliography{OCP-PINN}

\end{document}